\documentclass[onefignum,onetabnum]{siamart171218}

\usepackage{lipsum}
\usepackage{amsfonts}
\usepackage{graphicx}
\usepackage{epstopdf}
\usepackage{algorithmic}
\ifpdf
  \DeclareGraphicsExtensions{.eps,.pdf,.png,.jpg}
\else
  \DeclareGraphicsExtensions{.eps}
\fi

\newsiamremark{remark}{Remark}
\newsiamremark{hypothesis}{Hypothesis}
\crefname{hypothesis}{Hypothesis}{Hypotheses}
\newsiamthm{claim}{Claim}

\headers{One-dimensional ridge function approximation}{Andrew Glaws and Paul G. Constantine}

\title{A Lanczos-Stieltjes method for one-dimensional ridge function approximation and integration\thanks{\textbf{Funding:} This work was funded by the DARPA DSO Enabling Quantification of Uncertainty in Physical Systems program.}}

\author{Andrew Glaws\thanks{Department of Computer Science, University of Colorado Boulder,
  (\email{andrew.glaws@colorado.edu}).}
\and Paul G. Constantine\thanks{Department of Computer Science, University of Colorado Boulder,
  (\email{paul.constantine@colorado.edu}).}
}

\usepackage{amsopn}
\DeclareMathOperator{\diag}{diag}

\graphicspath{{./}{./figs/}}

  \def\clap#1{\hbox to 0pt{\hss#1\hss}}

\providecommand{\mat}[1]{\bm{#1}}%
\renewcommand{\vec}[1]{\mathbf{#1}}
\newcommand{\vecalt}[1]{\bm{#1}}

\providecommand{\mA}{\ensuremath{\mat{A}}}
\providecommand{\mB}{\ensuremath{\mat{B}}}
\providecommand{\mC}{\ensuremath{\mat{C}}}

\providecommand{\mI}{\ensuremath{\mat{I}}}
\providecommand{\mJ}{\ensuremath{\mat{J}}}

\providecommand{\mQ}{\ensuremath{\mat{Q}}}

\providecommand{\mT}{\ensuremath{\mat{T}}}

\providecommand{\mV}{\ensuremath{\mat{V}}}

\providecommand{\mLambda}{\ensuremath{\mat{\Lambda}}}

\providecommand{\mzero}{\ensuremath{\mat{0}}}

\providecommand{\va}{\ensuremath{\vec{a}}}

\providecommand{\ve}{\ensuremath{\vec{e}}}

\providecommand{\vp}{\ensuremath{\vec{p}}}
\providecommand{\vq}{\ensuremath{\vec{q}}}

\providecommand{\vu}{\ensuremath{\vec{u}}}
\providecommand{\vv}{\ensuremath{\vec{v}}}
\providecommand{\vw}{\ensuremath{\vec{w}}}
\providecommand{\vx}{\ensuremath{\vec{x}}}

\providecommand{\valpha}{\ensuremath{\vecalt{\alpha}}}

\providecommand{\vphi}{\ensuremath{\vecalt{\phi}}}
\providecommand{\vxi}{\ensuremath{\vecalt{\xi}}}

\newcommand{\sU}{\mathcal{U}}

\newcommand{\Exp}[1]{\mathbb{E}\left[#1\right]}

\newcommand{\Var}[1]{\operatorname{Var}\left[#1\right]}

\newcommand{\bmat}[1]{\begin{bmatrix}#1\end{bmatrix}}

\newcommand{\sign}[1]{\mathrm{sign}\left(#1\right)}

\usepackage{amsmath,mathtools}
\usepackage{tabulary,booktabs}
\usepackage{thmtools,thm-restate}
\usepackage{algorithm}
\usepackage{subfig}
\usepackage{bm}
\usepackage{color}
\usepackage{hyperref}

\ifpdf
\hypersetup{
  pdftitle={Gaussian quadrature and polynomial approximation for one-dimensional ridge functions},
  pdfauthor={Andrew Glaws and Paul G. Constantine}
}
\fi

\begin{document}

\maketitle

% REQUIRED
\begin{abstract}
Many of the input-parameter-to-output-quantity-of-interest maps that arise in computational science admit a surprising low-dimensional structure, where the outputs vary primarily along a handful of directions in the high-dimensional input space. This type of structure is well modeled by a ridge function, which is a composition of a low-dimensional linear transformation with a nonlinear function. If the goal is to compute statistics of the output---e.g., as in uncertainty quantification or robust design---then one should exploit this low-dimensional structure, when present, to accelerate computations. We develop Gaussian quadrature and the associated polynomial approximation for one-dimensional ridge functions. The key elements of our method are (i) approximating the univariate density of the given linear combination of inputs by repeated convolutions and (ii) a Lanczos-Stieltjes method for constructing orthogonal polynomials and Gaussian quadrature.
\end{abstract}

% REQUIRED
\begin{keywords}
Gaussian quadrature, orthogonal polynomials, ridge functions
\end{keywords}

% REQUIRED
\begin{AMS}
42C05, 33C50, 41A55
\end{AMS}

\sloppy 

\section{Numerical methods for high-dimensional integration}
\label{sec:intro}

High-dimensional integration is a common problem in scientific computing arising from, for example, the need to estimate expectations in uncertainty quantification~\cite{SmithUQ2013,SullivanUQ2015} and robust design~\cite{Allen2004}, where physics-based simulation models contain parametric uncertainty. However, integration suffers from the curse of dimensionality~\cite{Traub98,Donoho00}; loosely, the amount of information (e.g., integrand evaluations) needed to estimate the integral to within a fixed tolerance grows exponentially with dimension (i.e., the number of independent variables affecting the integrand). Monte Carlo~\cite{Owen2013} is a popular method for high-dimensional integration, since its $\mathcal{O}(N^{-1/2})$ convergence rate is independent of dimension. However, this rate is also very slow, which precludes the possibility of high accuracy; to obtain $k$ accurate digits, one must sample the integrand $\mathcal{O}(10^{2k})$ times. There are many extensions to simple Monte Carlo that produce relatively higher accuracy (i.e., \emph{variance reduction}) including stratified sampling, control variates, and multi-level methods~\cite{giles2015}. Quasi-Monte Carlo methods~\cite{caflisch1998} have a superior dimension-independent convergence rate of $\mathcal{O}(N^{-1})$ for certain classes of functions, though high accuracy is still elusive. Interpolatory sparse grid integration~\cite{bungartz2004} converges rapidly---comparable to one-dimensional interpolatory integration rules~\cite{Trefethen2008}---for smooth functions with small mixed partial derivatives, which makes high accuracy possible. Common sparse grid rules contain a \emph{level} parameter where the nodes of the level $\ell-1$ sparse grid are a subset of the level $\ell$ sparse grid nodes, which enables practical numerical convergence studies; such rules are called \emph{nested}. Although sparse grids are theoretically optimal in a precise sense~\cite{Novak1996}, the number of nodes needed to advance from level $\ell-1$ to $\ell$ is still too large for many scientific computing applications with expensive integrands.

An alternative strategy is to approximate the integrand with an easy-to-integrate approximation, such as a multivariate polynomial~\cite{hyman2014}. For any polynomial expressed in an orthogonal (with respect to the integration measure) polynomial basis, the integral of the polynomial is the coefficient associated with the constant term; this fact has contributed to the popularity of so-called \emph{polynomial chaos} methods in uncertainty quantification~\cite{LeMaitre2010}. From a Monte Carlo perspective, an approximation-based approach is comparable to control variates~\cite{Owen2013}. The integration error is tied to the integrand approximation error. Unfortunately, polynomial and related approximations are also dimensionally cursed; approximation-based approaches effectively trade one intractable high-dimensional problem for another. But the general idea of identifying exploitable structure in the integrand remains appealing and motivates our approach. 

\subsection{Ridge structure in scientific computing applications}
\label{sub:ridges}

Many functions that map input parameters to an output quantity of interest found in scientific computing applications admit a surprising one-dimensional structure we call \emph{ridge structure}. A \emph{ridge function}~\cite{Pinkus15} is a function $f:\mathbb{R}^m\rightarrow\mathbb{R}$ of the form
\begin{equation}
\label{eq:ridge0}
f(\vx) \;=\; g(\va^\top\vx),
\end{equation}
where $\va\in\mathbb{R}^m$ is a constant vector called the \emph{ridge direction} and $g:\mathbb{R}\rightarrow\mathbb{R}$ is the \emph{ridge profile}~\cite{Mayer15}. Without loss of generality, we can assume $\va$ has unit 2-norm. Ridge functions are constant along directions in the domain that are orthogonal to $\va$; for $\vu$ such that $\vu^\top\va=0$,
\begin{equation}
f(\vx+\vu) \;=\; g\left(\va^\top(\vx+\vu)\right) 
\;=\; g(\va^\top\vx) \;=\; f(\vx).
\end{equation}
Despite their name, ridge functions are not related to \emph{ridge regression}~\cite[Chapter 3.4]{Hastie2009}, which is a linear approximation scheme with 2-norm regularization on the coefficients. In contrast, ridge functions are basic objects in \emph{projection pursuit regression}~\cite{Friedman1981}, which approximates a function with a sum of ridge functions. There is a lot of literature on approximation by sums of ridge functions (see, e.g., chapters 22-24 of Cheney and Light~\cite{Cheney2000}), partly because one-layer feed-forward neural networks can be written as sums of ridge functions, where the ridge profile is commonly called the \emph{activation function}~\cite{Higham2018}.

There are several methods for assessing whether a given function has ridge structure; these methods start by estimating a ridge direction. The gradient of a differentiable ridge function points in the ridge direction, since $\nabla f(\vx)=\va\,g'(\va^\top\vx)$, where $g'$ is the derivative of the profile. If an application code includes algorithmic differentiation~\cite{Griewank2008} capabilities for sensitivity analysis or optimization, then a single simulation may reveal the ridge direction in a model's input/output map. When gradients are not available, one can use a set of function evaluations or point queries to estimate the ridge direction. For a well-behaved simulation, finite difference approximations of the partial derivatives may suffice to estimate a gradient and hence the ridge direction. Cohen et al.~\cite{Cohen12} analyze a greedy procedure for estimating the ridge direction $\va$ from point queries when all elements of $\va$ are non-negative and $\va$ is compressible. Fornasier et al.~\cite{Fornaiser12} extend the analysis to the case where $\va$ in \eqref{eq:ridge0} is a tall matrix with compressible columns; their method uses finite difference approximations of random directional derivatives. Tyagi and Cevher~\cite{Tyagi14} relax the compressibility assumption. These papers contain theoretical recovery guarantees and complexity estimates for their algorithms.

In computational science, Constantine et al.~\cite{Constantine2014} study the matrix $\mC=\mathbb{E}[\nabla f\,\nabla f^\top]$, whose dominant eigenspaces are \emph{active subspaces}, and noted that $\mC$ is rank-one if and only if $f$ is a continuously differentiable ridge function~\cite{Constantine17b}. Thus, when $\mC$ is rank-one, its first eigenvector is the ridge direction. In practice, if $\mC$ is nearly rank-one, then the associated $f$ may be well-approximated by a ridge function whose ridge direction is the first eigenvector of $\mC$. Numerical methods for estimating $\mC$ may use gradient capabilities when present~\cite{Constantine2015,Holodnak2018}, and several approximations have been developed when gradients are not available~\cite{Lewis2016,constantine2015computing,eftekhari2017learning}. Another approach for estimating a ridge direction involves fitting a ridge function model with a set of data comprised of point queries of the computational model's input/output map~\cite{Constantine17b,Hokanson18}. Related techniques use a Bayesian framework to derive posterior distributions on the ridge direction conditioned on data, and these approaches have been used in the related statistical subfield of \emph{computer experiments}~\cite{Gramacy12,Liu2017,Tripathy2016}.

In statistical regression, techniques for \emph{sufficient dimension reduction} (SDR)~\cite{Li2018} use given predictor/response pairs to find low-dimensional subspaces of the predictor (i.e., \emph{input}) space that are statistically sufficient to characterize the response (i.e., output). These techniques have been applied to computational science models for sensitivity analysis and dimension reduction~\cite{Cook94b,Li2016}, where the predictor/response pairs are point queries from the computational model and thus do not require gradients. Glaws, Constantine, and Cook showed how the associated inverse regression methods from SDR should be interpreted as estimating ridge directions in deterministic models~\cite{Glaws17a}. 

No matter how one estimates a ridge direction $\va$, one can assess the viability of the one-dimensional ridge function model by plotting $f(\vx)$ versus $\va^\top\vx$ using a set of point queries. If the plot reveals (near) functional structure---i.e., that the function $f$ appears to be well represented by a univariate function of $\va^\top\vx$---then the particular parameter-to-quantity-of-interest map has ridge structure. In statistical regression, such plots are called \emph{sufficient summary plots}~\cite{Cook98}. To avoid confusion with the precise notion of statistical sufficiency, which is not valid in the case where the data are derived from a deterministic function, we prefer to call these plots \emph{shadow plots}, since the name invokes the analogy of the surface's shadow along all but one direction.

\begin{figure}[!h]
\centering
\subfloat[Scramjet model~\cite{Constantine15a}]{
\includegraphics[width=0.31\textwidth]{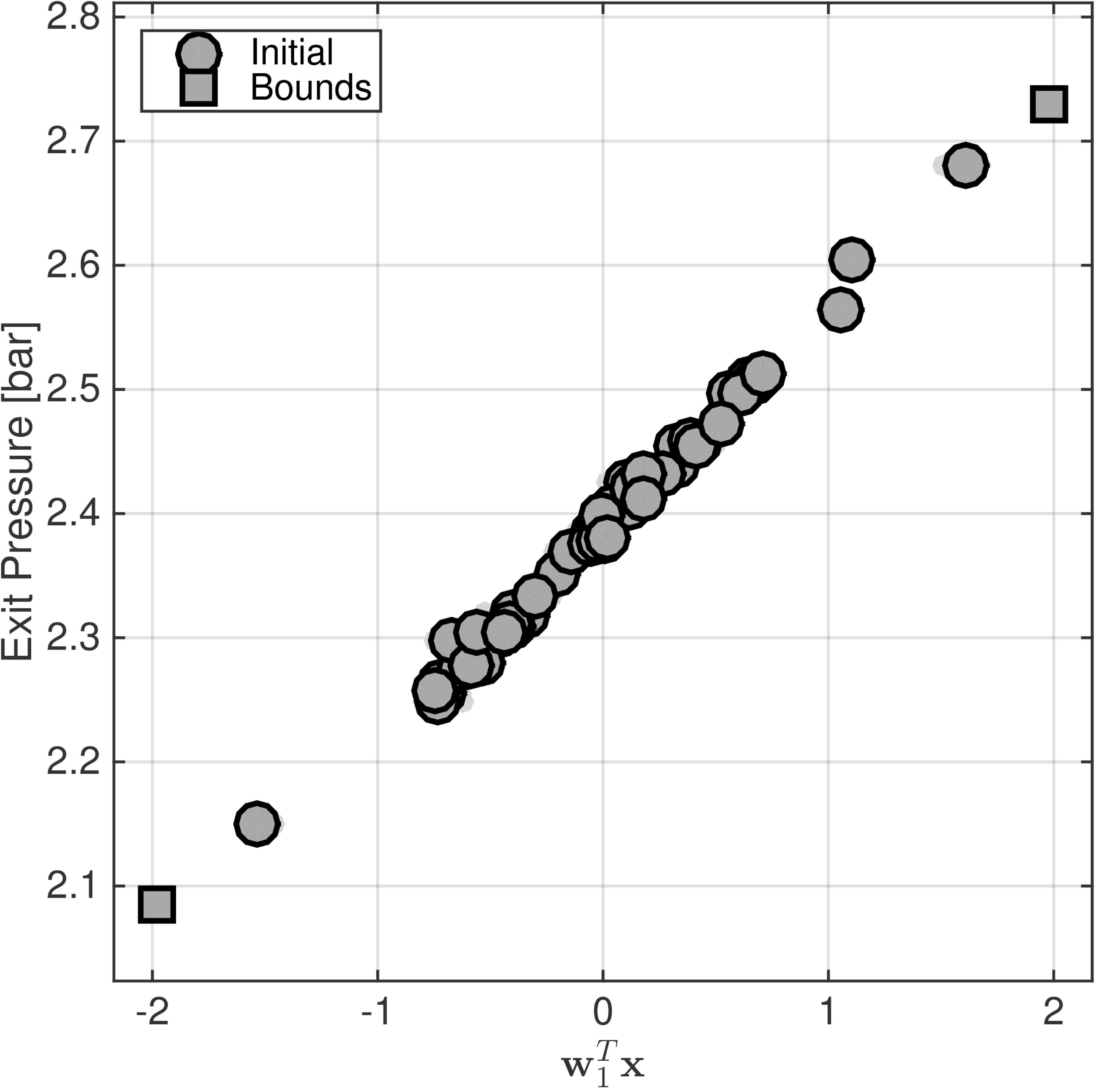}%
}
\hfil
\subfloat[Battery model~\cite{Constantine17a}]{
\includegraphics[width=0.31\textwidth]{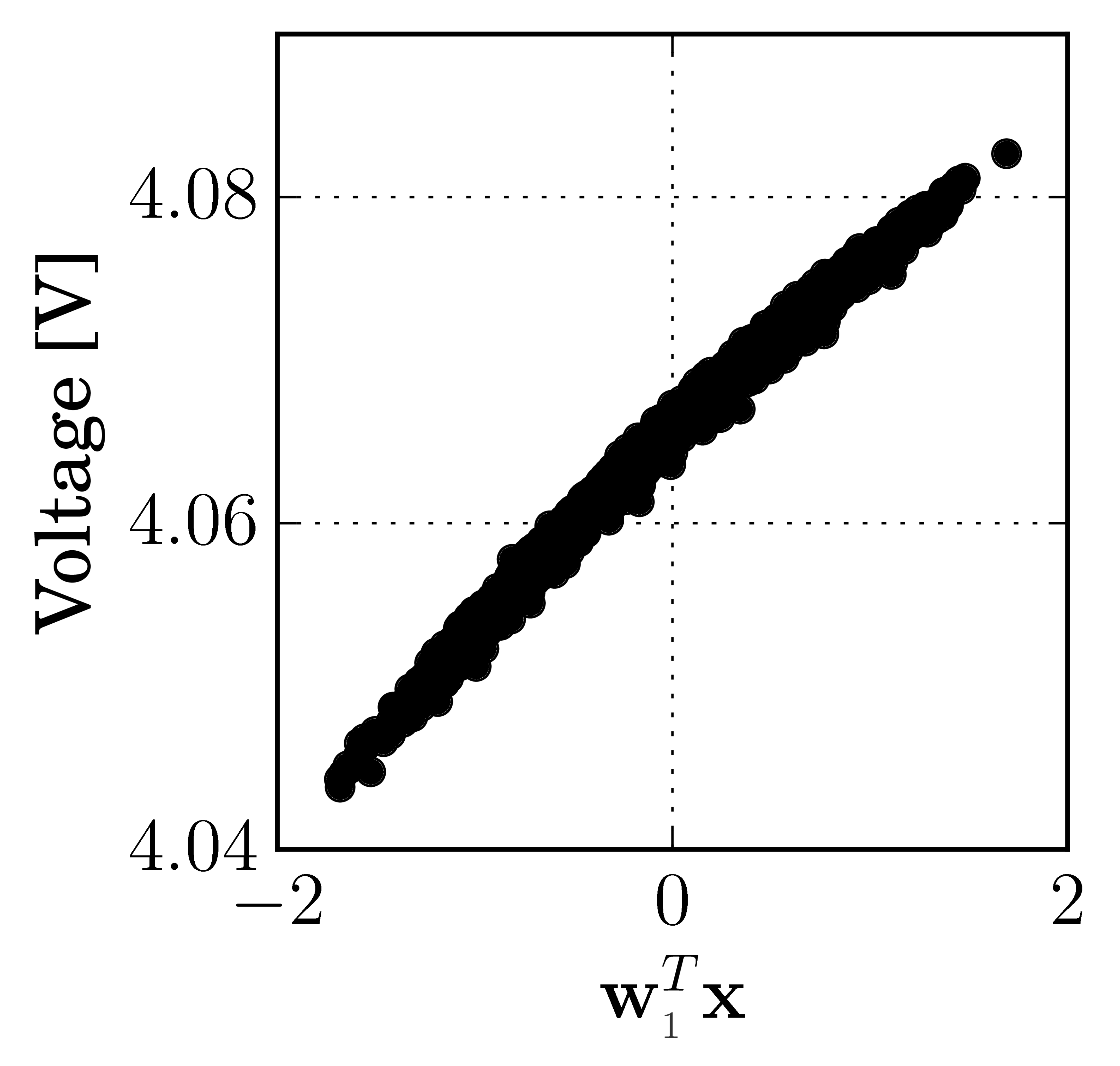}%
}
\hfil
\subfloat[Solar cell model~\cite{Constantine15b}]{
\includegraphics[width=0.31\textwidth]{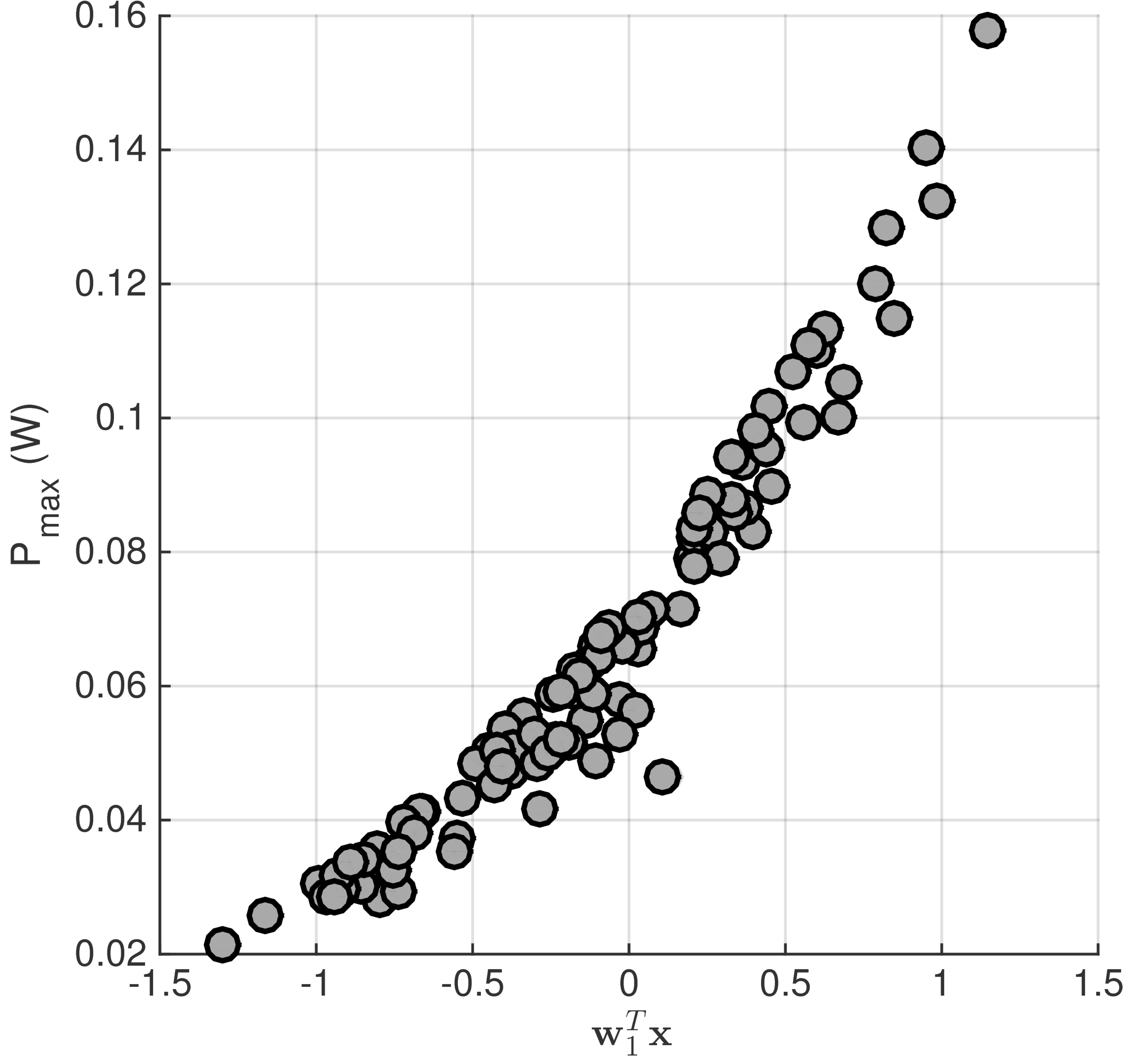}%
}\\
\subfloat[Reentry vehicle model~\cite{Cortesi2017}]{
\includegraphics[width=0.31\textwidth]{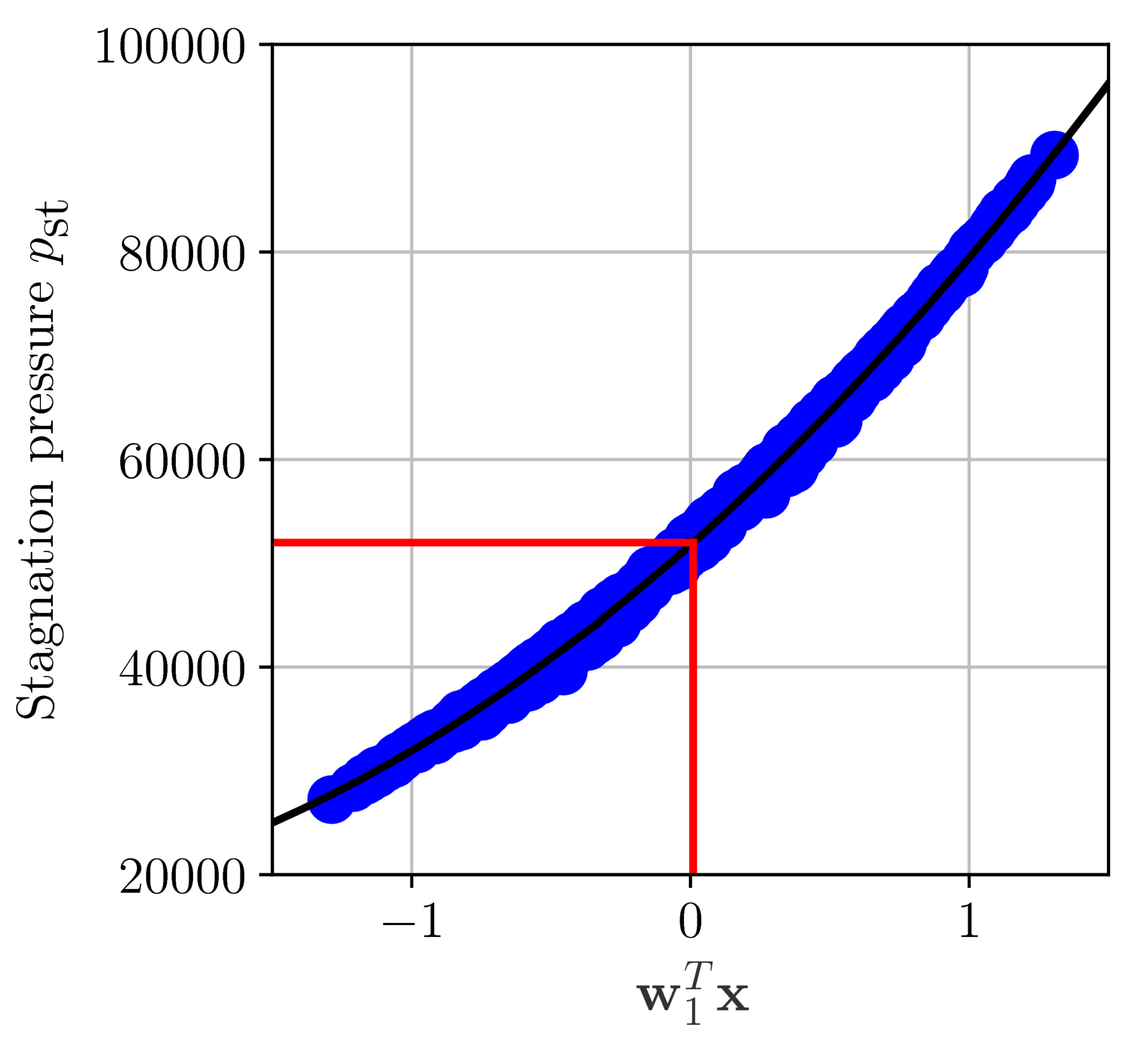}%
}
\hfil
\subfloat[Transonic airfoil model~\cite{Constantine15}]{
\includegraphics[width=0.31\textwidth]{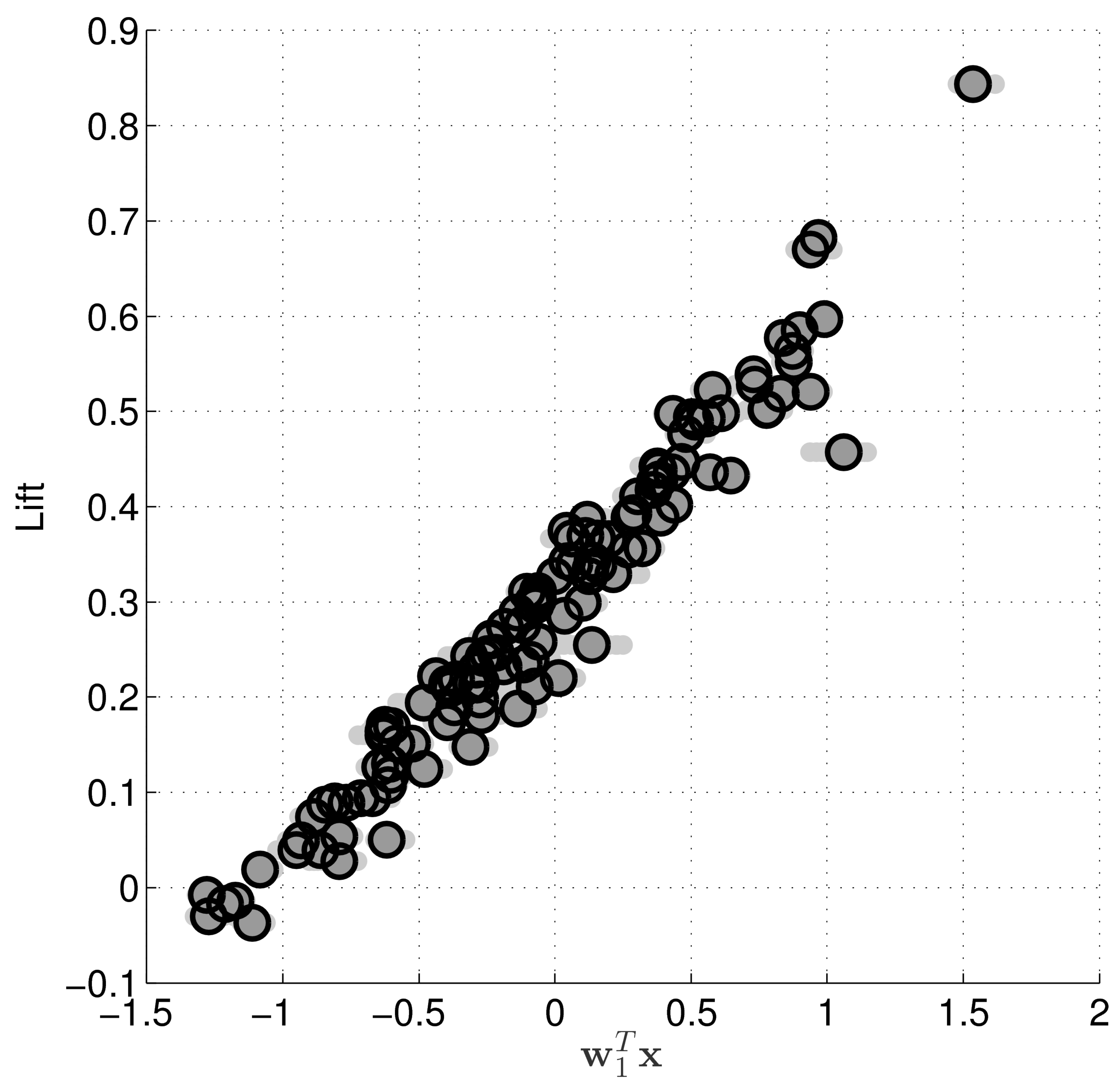}%
}
\hfil
\subfloat[MHD generator model~\cite{Glaws17b}]{
\includegraphics[width=0.31\textwidth]{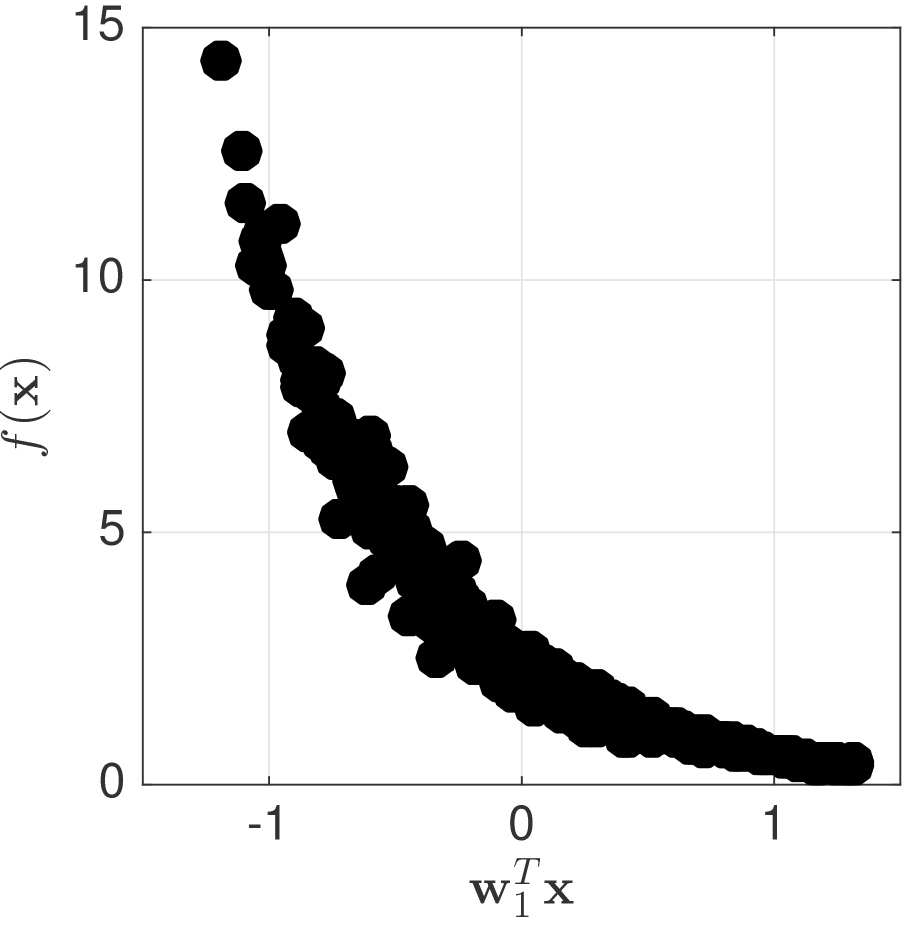}%
}
\caption{One-dimensional shadow plots of data derived from computational models across applications. The plots verify the near-one-dimensional structure in the input/output map.}
\label{fig:near_1d_ridge}
\end{figure}

Several of these exploratory approaches for estimating ridge directions and generating shadow plots have revealed one-dimensional ridge structure in a range of computational science applications, including wind farm modeling~\cite{King2018}, aerospace design~\cite{lukaczyk2014active,grey2018active,Berguin2014,Hu2016}, hypersonic vehicle modeling~\cite{Constantine15a}, turbine manufacturing~\cite{Seshadri2018,Beck2018}, hydrologic models~\cite{Jefferson16,Gilbert16}, and energy models~\cite{Glaws17b,Constantine17a,Constantine15b}. Figure \ref{fig:near_1d_ridge} shows shadow plots from several applications demonstrating the near one-dimensional ridge structure. To the best of our knowledge, there is neither a mathematical proof nor a first principles physical argument for why such structure is so prevalent across computational science applications. Nevertheless, the exploratory perspective continues to reveal such exploitable structure in real applications.

\subsection{Exploiting ridge structure for numerical integration}

The prevalence of one-dimensional ridge structure in applications justifies the question: if a function admits discoverable one-dimensional ridge structure, how can one exploit that structure to efficiently estimate the integral (e.g., the average) of the model's output quantity of interest? In other words, assume that a computational scientist can use one of the methods mentioned above to determine that a given model's input/output map is a ridge function---or is well-approximated by a ridge function---and the method estimates the ridge direction $\va$; how can the scientist exploit this information to efficiently estimate the map's integral?

Such integration requires approximation of the density induced by the linear combination $\va^\top \vx$. We compute this density efficiently using repeated convolutions. This approach is valid provided that the components of $\vx$ are independent, as is common in computational science applications. We then introduce a procedure for computing a generalized Gaussian quadrature rule for computing univariate integrals with respect to a known ridge direction in a high-dimensional space. The one-dimensional representation enables high accuracy with relatively few function evaluations. Since Gaussian quadrature is intimately tied to polynomial approximation, we obtain a univariate polynomial approximation of $f(\vx)$ using a basis of polynomials that are orthogonal with respect to the density of $\va^\top\vx$; this approximation can be used as a surrogate or response surface approximation of $f$. A related approach by Tsilifis uses the inverse empirical distribution to derive a polynomial approximation~\cite{Tsilifis2018}. However, this approach suffers from the relatively low accuracy of the empirical distribution; our convolution-based approach has no such limitation.  Additionally, low-dimensional Gaussian processes can be used with reduced-variance Monte Carlo integration~\cite{Gramacy12}, but this approach does not exhibit the same spectral convergence that is possible using the proposed Gaussian quadrature method.

The remainder of this paper proceeds as follows. Section \ref{sec:background} contains important background information, including discussions on ridge functions, Fourier expansions in terms of orthogonal polynomials, and the Lanczos and Stieltjes methods. In Section \ref{sec:exact_ridge}, we introduce the new Lanczos-Stieltjes method for one-dimensional ridge approximation and integration. Section \ref{sec:near_ridge} provides a heuristic extension of the method for approximate one-dimensional ridge functions, such as those shown in Figure \ref{fig:near_1d_ridge}.

%%%%

\section{Background}
\label{sec:background}

We represent the map from a computer model's physical inputs\footnote{We do not consider numerical parameters such as grid spacing or solver tolerances.} to its output quantity of interest as a scalar-valued function of $m$ independent variables,
\begin{equation}
\label{eq:y_fx}
y \;=\; f(\vx) , \qquad y \in \mathbb{R} , \quad \vx \in \mathbb{R}^m ,
\end{equation}
and we assume the input space is weighted by a given probability density function $p(\vx)$, which describes uncertainty in the model inputs. For simplicity, we assume this density is uniform over the $m$-dimensional hypercube $[ -1, 1]^m$ such that
\begin{equation}
\label{eq:uniform_dens}
p (\vx) \;=\;
\begin{cases} 
\frac{1}{2^m} & \text{if } ||\vx||_\infty \leq 1 , \\
0 & \text{otherwise.}
\end{cases}
\end{equation}
The independence is the important part; our approach extends easily to any product-type density function.

\subsection{Polynomial approximation}
\label{subsec:poly_approx}

Gaussian quadrature is intimately tied to approximation with orthogonal polynomials; for completeness, we provide a brief background. Assume that $f$ is square-integrable with respect to the input density $p(\vx)$. Then, $f$ admits a Fourier expansion in orthogonal polynomials with respect to $p$,
\begin{equation}
\label{eq:mD_exp}
y \;=\; f(\vx) \;=\; \sum_{|\valpha|=0}^{\infty} f_{\valpha} \, \psi_{\valpha} (\vx) ,
\end{equation}
where equality is denoted in the $L^2$ sense~\cite{Gautschi04}. The multivariate orthogonal polynomials $\psi_{\valpha} (\vx)$ are indexed by the multi-index $\valpha \in \mathbb{N}_0^m$ which denotes the degree of the polynomial with respect to each of the components of $\vx$. Since $p(\vx)$ is a product-type uniform density, the multivariate polynomials are
\begin{equation}
\psi_{\valpha} (\vx) \;=\; \prod_{i=1}^m \psi_{\alpha_i} (x_i) ,
\end{equation}
where each $\psi_{\alpha_i}$ is the univariate Legendre polynomial of degree $\alpha_i$. Without loss of generality, we also assume that the $\psi_{\valpha}$ are normalized so that the coefficients in \eqref{eq:mD_exp} are the inner product of $f$ with the appropriate polynomial,
\begin{equation}
\label{eq:mD_coeffs}
f_{\valpha} \;=\; \int f (\vx) \, \psi_{\valpha} (\vx) \, p(\vx) \, \text{d} \vx .
\end{equation}
This method of approximation by orthogonal polynomials also appears in the uncertainty quantification literature under the name polynomial chaos~\cite{Ghanem1991,Xiu02}.

In practice, we compute the \emph{pseudospectral expansion}~\cite{Constantine12b} of $f$ by truncating \eqref{eq:mD_exp} to include only polynomials of total degree $d$ or less and approximating the integral in \eqref{eq:mD_coeffs} numerically:
\begin{equation}
\label{eq:mD_pseudo_exp}
y \;=\; f(\vx) \;\approx\; \sum_{|\valpha| \leq d} \hat{f}_{\valpha} \, \psi_{\valpha} (\vx) , 
\quad \text{where} \quad
f_{\valpha} \;\approx\; \hat{f}_{\valpha} \;=\; \sum_{j=0}^{M-1} \omega_j \, f (\vxi_j) \, \psi_{\valpha} (\vxi_j) ,
\end{equation}
where $(\vxi_j, \omega_j)$, $j = 0, \dots, M-1$ denote the nodes and weights of an $M$-point numerical integration rule (e.g., tensor product Gauss quadrature) with respect to $p$.

Pseudospectral polynomial approximations serve as quick-to-evaluate surrogates for the original function. However, as the input dimension $m$ grows, the cost of constructing \eqref{eq:mD_pseudo_exp} quickly increases---a total degree $d$ polynomial in $m$ dimensions has ${m+d \choose d}$ coefficients. 

\subsection{Polynomial approximation for ridge functions}

\begin{figure*}[!ht]
\centering
\includegraphics[width=0.8\textwidth]{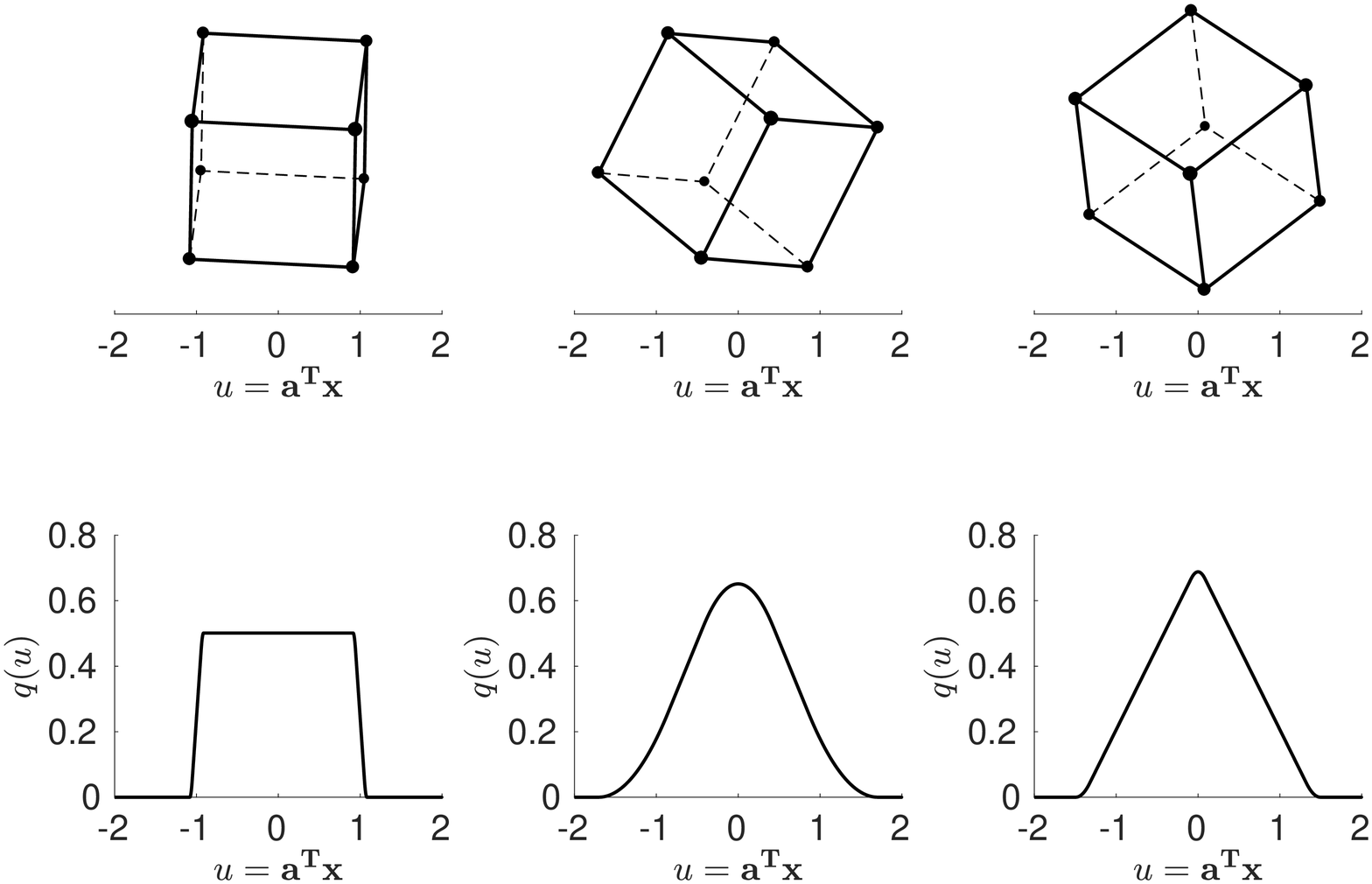}
\caption{Different density functions $q(u)$ induced by different vectors $\va \in \mathbb{R}^3 \setminus \{ \mzero \}$}
\label{fig:rotated_box}
\end{figure*}

To construct an orthogonal polynomial expansion of $g(\va^\top \vx)$ similar to \eqref{eq:mD_exp}, we must first understand the transformed input space. Let $u = \va^\top \vx$ denote the scalar-valued input of the ridge profile. The linear transform $\va^\top \vx$ induces a new density, which we denote $q(u)$. Figure \ref{fig:rotated_box} shows different rotations and projections, defined by different vectors $\va$, of the three-dimensional cube $[-1,1]^3$ and the resulting one-dimensional probability densities $q(u)$. 

We address the computation of $q(u)$ in Section \ref{sec:exact_ridge}. For now, assume $q(u)$ is known. We write the polynomial expansion of $g(u)$ as
\begin{equation}
\label{eq:1D_exp}
y \;=\; g(u) \;=\; \sum_{i = 0}^{\infty} g_i \, \phi_i (u) ,
\quad \text{where} \quad 
g_i \;=\; \int g(u) \, \phi_i (u) \, q(u) \, \text{d}u ,
\end{equation}
where the univariate polynomials $\phi_i$ are orthonormal with respect to $q$. By truncating the expansion and numerically approximating the coefficients, we obtain the pseudospectral approximation of the ridge profile,
\begin{equation}
\label{eq:1D_pseudo_exp}
y \;=\; g(u) \;\approx\; \sum_{i=0}^d \hat{g}_i \, \phi_i (u) ,
\quad \text{where} \quad 
\hat{g}_i \;=\; \sum_{j=0}^{M-1} \nu_j \, g(\lambda_j) \, \phi_i (\lambda_j) ,
\end{equation}
where $(\lambda_j \, , \, \nu_j)$ define a numerical integration rule with respect to $q$. %It simplifies later discussion if we express \eqref{eq:1D_pseudo_exp} in matrix-vector form as
%\begin{equation}
%y \;=\; g(u) \;=\; \hat{\vg}^\top \vphi(u) ,
%\quad \text{where} \quad
%\hat{\vg} \;=\; \mQ_u \, \mW_u \, \vg
%\end{equation}
%where
%\begin{equation}
%\begin{aligned}
%\mQ_u \;&=\; \bmat{\sqrt{\nu_0} \, \phi_0 (\lambda_0) & \dots & \sqrt{\nu_{N-1}} \, \phi_0 (\lambda_{N-1}) \\ \vdots & \ddots & \vdots \\ \sqrt{\nu_0} \, \phi_d (\lambda_0) & \dots & \sqrt{\nu_{N-1}} \, \phi_d (\lambda_{N-1})}, \\[1mm]
%\mW_u \;&=\; \bmat{\sqrt{\nu_0} & & \\ & \ddots & \\ & & \sqrt{\nu_{N-1}}},
%\;
%\vphi(u) \;=\; \bmat{\phi_0 (u) \\ \vdots \\ \phi_d (u)},
%\; \text{and} \;
%\vg \;=\; \bmat{g(\lambda_0) \\ \vdots \\ g(\lambda_{N-1})} .
%\end{aligned}
%\end{equation}

\begin{figure*}[!ht]
\centering
\includegraphics[width=0.75\textwidth]{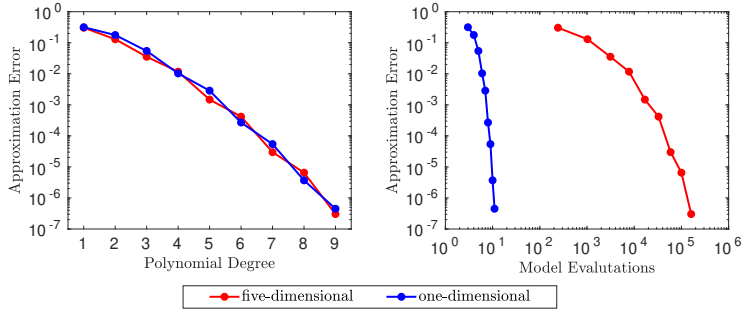}
\caption{A comparison of the cost and accuracy of the five-dimensional (in red) and the one-dimensional (in blue) polynomial approximations for a given function. They perform equally well in terms of polynomial degree, but the one-dimensional approximation uses four orders of magnitude fewer function evaluations.}
\label{fig:mD_vs_1D}
\end{figure*}

Constructing the pseudospectral polynomial expansion of the one-dimensional ridge profile $g(u)$ significantly reduces the number of function evaluations required compared to working in the $m$-dimensional space. Figure \ref{fig:mD_vs_1D} shows an example of approximating a smooth ridge function of five variables using either a multivariate polynomial in all five variables or a univariate polynomial along the ridge direction; unsurprisingly, the univariate approximation uses four orders of magnitude fewer function evaluations. However, \eqref{eq:1D_pseudo_exp} requires knowledge of the orthonormal polynomials $\phi_i$ and an integration rule $(\lambda_j \, , \, \nu_j)$ with respect to $q$, as well as the ability to evaluate the ridge profile, which may not be known in closed form, at the $\lambda_j$'s. We address this latter issue in Section \ref{sec:exact_ridge}. In the next section, we discuss how to obtain the orthonormal polynomials and integration nodes/weights using the Stieltjes and Lanczos iterative methods. 

We next review the Stieltjes and Lanczos iterative methods. We discuss the key components of each and show that, under certain conditions, the Lanczos method is a discrete approximation to the Stieltjes procedure. These methods and the relationships between them have been rigorously studied~\cite{Liesen13, Gautschi04, Golub10}; the discussion in the following section is based on these references.

\subsection{The Stieltjes procedure}
\label{subsec:stieltjes_lanczos}

The Stieltjes procedure---given in Algorithm \ref{alg:Stieltjes}---is a method for iteratively constructing a sequence of orthonormal polynomials $\{ \phi_0 (u), \phi_1 (u), \phi_2 (u), \dots \}$ with respect to a given density~\cite{Stieltjes84}. Step (iv) in Algorithm \ref{alg:Stieltjes} contains the three-term recurrence relationship that must be satisfied by any sequence of orthonormal polynomials. 

\begin{algorithm}
\caption{Stieltjes procedure~\cite[Section 2.2.3.1]{Gautschi04}} \label{alg:Stieltjes}
\textbf{Given:} probability density function $q(u)$ \\
\textbf{Assumptions:} $\phi_{-1} (u) = 0$ and $\tilde{\phi}_0 (u) = 1$ \\
For $i \;=\; 0, 1, 2, \dots$
\begin{enumerate}
\item[(i)] $\displaystyle \beta_i \;=\; \int \tilde{\phi}_i (u)^2 \, q(u) \, \text{d}u$
\item[(ii)] $\displaystyle \phi_i (u) \;=\; \tilde{\phi}_i (u) \, / \, \beta_i$
\item[(iii)] $\displaystyle \alpha_i \;=\; \int u \, \phi_i (u)^2 \, q(u) \, \text{d}u$
\item[(iv)] $\displaystyle \tilde{\phi}_{i+1} (u) \;=\; (u - \alpha_i) \, \phi_i (u) - \beta_i \, \phi_{i-1}(u)$
\end{enumerate}
\textbf{Output:} the orthonormal polynomials $\{ \phi_0 (u), \phi_1 (u), \phi_2 (u), \dots \}$ and recurrence coefficients $\alpha_i$, $\beta_i$ for $i = 0,1,2,\dots$
\end{algorithm}

By stopping the Stieltjes procedure after $d+1$ terms, we can rewrite the three-term recurrence relationship in vector form as
\begin{equation}
\label{eq:three_term_vec}
u \, \vphi (u) \;=\; \mJ \, \vphi (u) + \beta_{d+1} \, \phi_{d+1} (u) \ve_{d+1} ,
\end{equation}
where $\vphi (u) = \left[ \, \phi_0 (u) \, , \, \phi_1 (u) \, , \, \dots \, , \, \phi_d (u) \, \right]^\top$, $\ve_{d+1} \in \mathbb{R}^{d+1}$ is a vector of zeros with a one in the last entry, and the matrix $\mJ \in \mathbb{R}^{(d+1) \times (d+1)}$---referred to as the \emph{Jacobi matrix}---is a symmetric, tridiagonal matrix of recurrence coefficients,
\begin{equation}
\label{eq:J_Jacobi}
\mJ \;=\; \bmat{
\alpha_0 & \beta_1 &  &  &  \\
\beta_1 & \alpha_1 & \beta_2 &  &  \\
 & \ddots & \ddots & \ddots &  \\
 &  & \beta_{d-1} & \alpha_{d-1} & \beta_d \\
 &  &  & \beta_d & \alpha_d } .
\end{equation}
Let $\mJ \;=\; \mQ \, \mLambda \, \mQ^\top$ denote the eigendecomposition of $\mJ$. From \eqref{eq:three_term_vec}, the $d+1$ eigenvalues of $\mJ$ are the zeros of the ($d+1$)-degree orthonormal polynomial $\phi_{d+1}$. The normalized eigenvector associated with the eigenvalue $\lambda_j$ has the form
\begin{equation}
\left( \mQ \right)_j \;=\; \frac{\vphi (\lambda_j)}{\sqrt{\vphi (\lambda_j)^\top \vphi (\lambda_j)}} ,
\end{equation}
where $(\cdot)_j$ denotes the $j$th column of the given matrix.

\subsubsection{Gaussian quadrature from the Stieltjes procedure}
\label{sub:gq_Stieltjes}

The $(d+1)$-point Gaussian quadrature rule with respect to $q$ can be obtained from the eigendecomposition of $\mJ$ from \eqref{eq:J_Jacobi}~\cite{Golub69}. The quadrature nodes are the eigenvalues of $\mJ$ and the associated quadrature weights are the square of the first entry of the associated normalized eigenvector,
\begin{equation}
\nu_j 
\;=\; \left( \mQ \right)_{0, j}^2
\;=\; \frac{1}{\vphi (\lambda_j)^\top \vphi (\lambda_j)} .
\end{equation}

\subsection{The Lanczos procedure}
\label{sub:Lanczos}

The Lanczos algorithm was introduced as an iterative approach to estimating eigenvalues of large symmetric matrices~\cite{lanczos50}. Given symmetric $\mA \in \mathbb{R}^{N \times N}$, the Lanczos algorithm constructs the system
\begin{equation}
\mA \, \mV = \mV \, \mT + \beta_{d+1} \, \vv_{d+1}^\top \, \ve_{d+1} ,
\end{equation}
where $\mT \in \mathbb{R}^{(d+1) \times (d+1)}$ is a symmetric, tridiagonal matrix of recurrence coefficients, $\mV \in \mathbb{R}^{N \times (d+1)}$ contains the Lanczos vectors, and $\ve_{d+1} \in \mathbb{R}^{d+1}$ is a vector of zeros with a one in the last entry. The eigenvalues of $\mT$ approximate those of $\mA$, and $\mV$ transforms the eigenvectors of $\mT$ into approximate eigenvectors of $\mA$. We consider conditions under which the Lanczos algorithm serves as a discrete approximation to the Stieltjes procedure. Algorithm \ref{alg:Lanczos} contains the Lanczos algorithm. For notational convenience, we put the outputs of Algorithm \ref{alg:Lanczos} into matrices. Define the matrix $\mV$ of Lanczos vectors as
\begin{equation}
\label{eq:V_lanczos}
\mV
\;=\;
\bmat{
\vline & \vline & & \vline \\
\vv_0 & \vv_1 & \dots & \vv_{d-1} \\
\vline & \vline & & \vline },
\end{equation}
and the symmetric, tridiagonal Jacobi matrix $\mT$ as
\begin{equation}
\label{eq:T_Jacobi}
\mT
\;=\;
\bmat{
\alpha_0 & \beta_1 & & & \\
\beta_1 & \alpha_1 & \beta_2 & & \\
 & \ddots & \ddots & \ddots & \\
 & & \beta_{d-2} & \alpha_{d-2} & \beta_{d-1} \\
 & & & \beta_{d-1} & \alpha_{d-1}} .
\end{equation}

\begin{algorithm}
\caption{Lanczos algorithm~\cite[Section 3.1.7.1]{Gautschi04}} \label{alg:Lanczos}
\textbf{Given:} an $N \times N$ symmetric matrix $\mA$ \\
\textbf{Assumptions:} $\vv_{-1} = \mzero \in \mathbb{R}^N$ and $\tilde{\vv}_0 \in \mathbb{R}^N \setminus \{ \mzero \}$

For $i \;=\; 0,1\dots,d-1$,
\begin{enumerate}
\item $\beta_i \;=\; \sqrt{\tilde{\vv}_i^\top \tilde{\vv}_i}$
\item $\vv_i \;=\; \tilde{\vv}_i \, / \, \beta_i$
\item $\alpha_i \;=\; \vv_i^\top \mA \, \vv_i$
\item $\tilde{\vv}_{i+1} \;=\; (\mA - \alpha_i \, \mI) \, \vv_i - \beta_{i-1} \, \vv_{i-1}$
\end{enumerate}

\end{algorithm}

\subsubsection{Gaussian quadrature from Lanczos}
\label{sub:gq_Lanczos}

Let $(u_j, \upsilon_j)$, $j = 0, \dots, N-1$ be the nodes and weights of an $N$-point numerical integration rule with respect to a given density function $q(u)$. These nodes and weights define a discrete approximation of $q(u)$, which we denote by $q^{(N)} (u)$. Performing Algorithm \ref{alg:Lanczos} on the diagonal matrix
\begin{equation}
\mA = \bmat{u_0 & & \\ & \ddots & \\ & & u_{N-1}}
\end{equation}
with starting vector
\begin{equation}
\tilde{\vv}_0 = \bmat{\sqrt{\upsilon_0} \\ \vdots \\ \sqrt{\upsilon_{N-1}}}
\end{equation}
is equivalent to performing Algorithm \ref{alg:Stieltjes} on the discrete density function $q^{(N)} (u)$~\cite{Constantine2012}. The recurrence coefficients in $\mT$ from \eqref{eq:T_Jacobi} converge to those in $\mJ$ from \eqref{eq:J_Jacobi} as $N$ goes to infinity~\cite{Gautschi04}. Therefore, we can use these recurrence coefficients to produce approximations to the orthonormal polynomials $\{ \phi_0 (u), \phi_1 (u), \phi_2 (u), \dots\}$. Additionally, the eigendecomposition of $\mT$ provides us with an approximate Gaussian quadrature rule with respect to $q(u)$.

\section{Integration and approximation for ridge functions}
\label{sec:exact_ridge}

Assume that $f$ is a ridge function,
\begin{equation}
\label{eq:y_fx_gaTx}
y \;=\; f(\vx) \;=\; g(u), 
\quad \text{where} \quad
u \;=\; \va^\top \vx .
\end{equation} We want to build a Gaussian quadrature rule and polynomial approximation with respect to the induced density function $q(u)$. 

%Recall from \eqref{eq:y_fx_gAx} that we assume $f(\vx)$ is a one-dimensional ridge function,
%\begin{equation}
%\label{eq:y_fx_gAx2}
%y \;=\; f(\vx) \;=\; g(u) ,
%\end{equation}
%where $u = \va^\top \vx$. We want to build a pseudospectral expansion of $g$ 
%\begin{equation}
%\label{eq:1D_pseudo_exp2}
%y \;=\; g(u) \;\approx\; \sum_{i=0}^d \hat{g}_i \, \phi_i (u) ,
%\quad
%\hat{g}_i \;=\; \sum_{j=0}^{M-1} \nu_j \, g(\lambda_j) \, \phi_i (\lambda_j) ,
%\end{equation}
%where $\phi_i$ denote the orthonormal polynomials with respect to the induced density function $q(u)$.

\subsection{Computing the density of a linear combination of independent variables}

We compute $q(u)$ using a \emph{convolution of probability densities}~\cite[Ch. 4]{Billingsley86}. Consider two independent random variables $x_0$ and $x_1$ with density functions $p_0$ and $p_1$, respectively. The density function of $u = x_0 + x_1$ is 
\begin{equation}
\label{eq:convolve_dens}
q(u) \;=\; (p_0 * p_1) (u) \;=\; \int p_0 (t) \, p_1 (u - t) \, \text{d}t .
\end{equation}
Equation \eqref{eq:convolve_dens} is the convolution of $p_0$ and $p_1$.

Recall from \eqref{eq:uniform_dens} that we assume the input space is weighted by a uniform density over the hypercube, $\vx \sim \sU ([-1, 1]^m)$. By independence, we have that $p(\vx) = p_0(x_0) \, p_1(x_1) \, \dots \, p_{m-1} (x_{m-1})$ with each
\begin{equation}
\label{eq:uniform_dens_1D}
p_i (x_i) \;=\;
\begin{cases} 
\frac{1}{2} & \text{if } |x_i| \leq 1 , \\
0 & \text{otherwise.}
\end{cases}
\end{equation}
We can write the linear transform $u = \va^\top \vx = a_0 \, x_0 + \dots + a_{m-1} x_{m-1}$ and recognize that $a_i \, x_i$ is uniformly distributed in the interval $[-a_i,a_i]$. Thus, we can obtain $q(u)$ by iteratively applying convolutions to each $a_i x_i$. In practice, we approximate the integral in \eqref{eq:convolve_dens} using a trapezoidal rule since the integrands in \eqref{eq:convolve_dens} are not smooth. However, this integral is one-dimensional and does not require us to evaluate the computational model ($f(\vx)$ from \eqref{eq:y_fx}), so we can use a high density of points to approximate the convolution. 

Algorithm \ref{alg:convolve_dens} details the process of approximating $q(u)$ using iterative convolutions. Step (1) finds the range of $u$. Equation \eqref{eq:ubounds} is justified by noting that the maximizer of $\va^\top\vx$ over the hypercube is $\sign{\va}$, and the minimizer of $\va^\top\vx$ is $-\sign{\va}$. If $a_i = 0$ for some $i$, then the corresponding $x_i$ has no influence on the model output and we can remove this variable. 
%Also note that we define $N$ equally-spaced points along the one-dimensional interval, where $N$ is odd. This requirement is an artifact of the numerical approximation of \eqref{eq:convolve_dens}. 

\begin{algorithm}[!ht]
\caption{Discrete convolution of densities} \label{alg:convolve_dens}
\textbf{Given:} the vector $\va \in \mathbb{R}^m \setminus \{ \mzero \}$ \\
\textbf{Assumptions:} $\vx \sim \sU ([-1, 1]^m)$
\begin{enumerate}
\item Find the inputs of the one-dimensional interval
\begin{equation}
\label{eq:ubounds}
u_\ell \;=\; \va^\top \sign{-\va} , \quad
u_r \;=\; \va^\top \sign{\va} , 
\end{equation}
and define $N$ (where $N$ is odd) equally-spaced points along the interval
\begin{equation}
u_j \;=\; u_\ell + j \, \Delta u, 
\quad j \;=\; 0,\dots,N-1 ,
\end{equation}
where $\Delta u = (u_r - u_\ell)/(N-1)$.
\item Initialize the vector $\vq = \bmat{q(u_0) & \dots & q(u_{N-1})}^\top$ where
\begin{equation}
q(u_j) \;=\; \begin{cases} 
1/(2 \, a_0) & \text{if } |u_j| \leq a_0 , \\
0 & \text{otherwise.}
\end{cases}
\end{equation}
\item For $i = 1, \dots, m-1$
\begin{enumerate}
\item[(i)] Define $\vp = \bmat{p(u_0) & \dots & p(u_{N-1})}^\top$ where\begin{equation}
p(u_j) \;=\; \begin{cases} 
1/(2 \, a_i) & \text{if } |u_j| \leq a_i , \\
0 & \text{otherwise.}
\end{cases}
\end{equation}
\item[(ii)] For $j = 0, \dots, N-1$, define 
\begin{equation}
\begin{aligned}
k_0 \;&=\; \max \left\{ 0, j - \frac{N-1}{2} \right\} , \\
k_1 \;&=\; \min \left\{ \frac{N-1}{2} + j, N-1 \right\} ,
\end{aligned}
\end{equation}
and compute
\begin{equation}
\label{eq:discrete_conv}
q_j \;=\; \sum_{k=k_0}^{k_1} q_k \, p_{\frac{N-1}{2} - j + k}  
\end{equation}
\end{enumerate}
\end{enumerate}
NOTE: skip any $i$ for which $a_i = 0$ \\
\textbf{Output:} $\vq = \bmat{q_0 & \dots & q_{N-1}}^\top$
\end{algorithm}

Each element $q_i$ of the output vector $\vq$ from Algorithm \ref{alg:convolve_dens} approximates the density $q$ at the point $u_i$. The sum in \eqref{eq:discrete_conv} is a discrete approximation to the convolution \eqref{eq:convolve_dens}. Modern implementations of discrete convolution, such as Matlab's \texttt{conv}, use the fast Fourier transform that enables scaling to very large $N$, which in turn controls the error in the $q_i$'s. We stress that no computational model evaluations (i.e., evaluations of $f$ from \eqref{eq:y_fx}) are needed to compute the $q_i$'s. 

\subsection{Obtaining the Gaussian quadrature for the convolved density}

To obtain the Gaussian quadrature rule for $q(u)$, we apply the Lanczos procedure (Section \ref{sub:gq_Lanczos}) to the diagonal matrix $\diag([u_0,\dots,u_{N-1}])$ with starting vector $[q_0^{1/2},\dots,q_{N-1}^{1/2}]^\top$ computed with  Algorithm \ref{alg:convolve_dens}. Effectively, this runs the Stieltjes' procedure (Section \ref{sub:gq_Stieltjes}) using a discrete inner product defined by the trapezoidal rule with $N$ points. Since the matrix is diagonal, all computations use $\mathcal{O}(N)$ operations; thus, we can use very large $N$ to control errors in the quadrature rule calculation. Again, we stress that no model evaluations are needed. 

\subsection{Evaluating the function at the Gaussian quadrature points} 

\begin{figure*}[!ht]
\centering
\includegraphics[width=0.35\textwidth]{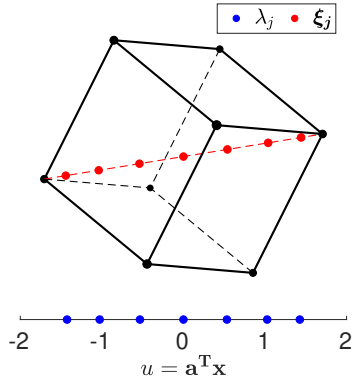}
\caption{The projection of the one-dimensional quadrature nodes $\lambda_j$ into $m$-dimensional space.}
\label{fig:line_projection}
\end{figure*}

We next address the issue of evaluating the ridge profile at the Gaussian quadrature nodes, $g(\lambda_j)$. Consider the linear transformation $\vxi_j = \lambda_j \, \va$, and recognize that 
\begin{equation}
\label{eq:1dQuad_to_mdQuad}
f(\vxi_j) \;=\; f(\lambda_j \, \va) \;=\; g(\va^\top (\lambda_j \, \va)) \;=\; g(\lambda_j \, \va^\top \va) \;=\; g(\lambda_j),
\end{equation}
since $\va$ is assumed to be normalized. Thus, we can evaluate $f(\vxi_j)$ in place of $g(\lambda_j)$, provided that $\vxi_j \in [-1,1]^m$. To ensure we find a valid point at which to evaluate $f$, consider the projection of the $m$-dimensional hypercube down to a one-dimensional domain via $u = \va^\top \vx$. In general, the endpoints of this one-dimensional interval are defined by two opposing corners of the hypercube. The endpoints of the interval are
\begin{equation}
\label{eq:u_corners}
u_\ell \;=\; \va^\top \sign{-\va}
\quad \text{and} \quad
u_r \;=\; \va^\top \sign{\va} ,
\end{equation}
and the corresponding corners of the hypercube are
\begin{equation}
\vx_\ell \;=\; \sign{-\va}
\quad \text{and} \quad
\vx_r \;=\; \sign{\va}.
\end{equation}
A line in $m$-dimensional space that connects these two corners of the hypercube is guaranteed to be contained within the hypercube. By projecting the one-dimensional quadrature points $\lambda_j$ onto that line, we ensure that $\vxi_j \in [-1, 1]^m$. 
We do this by
\begin{equation}
\label{eq:xi_j}
\vxi_j \;=\; ( 1 - \gamma_j ) \, \vx_\ell + \gamma_j \, \vx_r ,
\end{equation}
where $\gamma_j = (\lambda_j - u_\ell) / (u_r - u_\ell)$. Figure \ref{fig:line_projection} illustrates this process for a three-dimensional cube.

Algorithm \ref{alg:ridge_approx} summarizes the steps for building a pseudospectral polynomial approximation of a ridge function. The approximate integral of $f$ is the first pseudospectral coefficient:
\begin{equation}
\hat{g}_0 \;=\; \sum_{j=0}^{d} \nu_j \, g(\lambda_j) \;\approx\; \int g(u) \, q(u) \, \text{d}u \;=\; \int f(\vx) \, p(\vx) \, \text{d} \vx .
\end{equation}

\begin{algorithm}[!ht] 
\caption{Univariate polynomial approximation of a ridge function}
\label{alg:ridge_approx}
\textbf{Given:} function $f: \mathbb{R}^m \rightarrow \mathbb{R}$ and unit vector $\va \in \mathbb{R}^m \setminus \{ \mzero \}$ such that
$$y \;=\; f(\vx) \;=\; g(\va^\top \vx)$$
for some unknown $g: \mathbb{R} \rightarrow \mathbb{R}$ \\
\textbf{Assumptions:} $\vx \sim \sU ([-1, 1]^m)$ and $f$ is square-integrable with respect to the input density $p(\vx)$
\begin{enumerate}
\item Perform Algorithm \ref{alg:convolve_dens} to obtain $\vq = \bmat{q(u_0) & \dots & q(u_{N-1})}^\top$.
\item Perform Algorithm \ref{alg:Lanczos} on
\begin{equation}
\mA \;=\; \bmat{u_0 & & \\ & \ddots & \\ & & u_{N-1}} ,
\quad
\tilde{\vv}_0 \;=\; \bmat{\sqrt{q(u_0)} \\ \vdots \\ \sqrt{q(u_{N-1})}} ,
\end{equation}
to obtain the Jacobi matrix $\mT$.
\item Take the eigendecomposition of $\mT$,
\begin{equation}
\mT \;=\; \mQ \, \mLambda \, \mQ^\top ,
\end{equation}
where the eigenvalues are $\{ \lambda_0, \lambda_1, \dots, \lambda_d \}$ and the eigenvectors are normalized. Define $\nu_j = (\mQ)_{0, j}^2$ for $j = 0,\dots,d$.
\item For $j = 0,\dots,d$, compute
\begin{equation}
\vxi_j \;=\; (1 - \gamma_j) \, \vx_\ell + \gamma_j \, \vx_r ,
\quad \text{for} \quad 
\gamma_j = (\lambda_j - u_\ell) \, / \, (u_r - u_\ell) ,
\end{equation}
where $\vx_\ell = \sign{-\va}$, $\vx_r = \sign{\va}$, $u_\ell = \va^\top \sign{-\va}$, and $u_r = \va^\top \sign{\va}$.
\item Compute the pseudospectral coefficients 
\begin{equation}
\hat{g}_i \;=\; \sum_{j=0}^{d} \nu_j \, f(\vxi_j) \, \phi_i (\lambda_j) .
\end{equation}
\item Build the pseudospectral expansion 
\begin{equation}
y \;=\; g(u) \;\approx\; \sum_{i=0}^d \hat{g}_i \, \phi_i (u) .
\end{equation}
\end{enumerate}
\end{algorithm}

Algorithm \ref{alg:ridge_approx} contains two levels of approximation: (i) the discrete approximation of $q(u)$ using an $N$-point trapezoidal rule and (ii) the number $d+1$ of Lanczos iterations performed. The number of Lanczos iterations corresponds to the number of terms in the polynomial approximation of $g$ and the number of Gaussian quadrature nodes. The latter is important as this is the number of model evaluations required to compute the pseudospectral coefficients and the estimated integral. In general, we should choose $N \gg d$ since the approximation of $q$ at the $N$ trapezoidal points does not require any function evaluations, which are typically the most expensive step. In the next section, we numerically study this two-level approximation.

\subsection{Numerical study}
\label{subsec:exact_ridge_numerics}

In this section, we numerically study the behavior of Algorithm \ref{alg:ridge_approx} for
\begin{equation}
\label{eq:ex_exactRidge}
y \;=\; \sin \left( 2 \, \pi \, (\va^\top \vx) \right) + \cos \left( \frac{\pi}{2} (\va^\top \vx) \right) ,
\qquad
\vx \in \mathbb{R}^{25} .
\end{equation}
We assume $\vx \sim \sU ([-1, 1]^{25})$. Notice that \eqref{eq:ex_exactRidge} is an exact one-dimensional ridge function.

\begin{figure*}[!ht]
\centering
\includegraphics[width=1\textwidth]{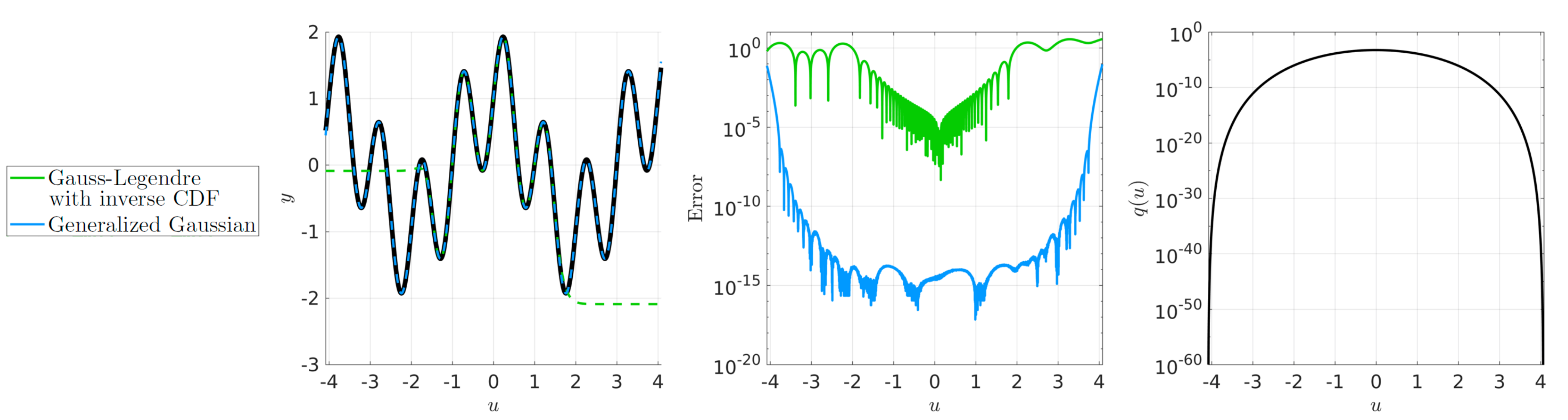}
\caption{The results of Algorithm \ref{alg:ridge_approx} applied to \eqref{eq:ex_exactRidge}. The first plot shows the true ridge profile (in black) and the $d = 50$ pseudospectral polynomial approximation (in blue) and the Gauss-Legendre approximation from Tsilifis~\cite{Tsilifis2018} (in green). The second plot contains the absolute error of the approximations, and the third plot shows the density $q(u)$.}
\label{fig:eyeball_norm}
\end{figure*}

Figure \ref{fig:eyeball_norm} contains the results from Algorithm \ref{alg:ridge_approx} for $N = 10,000$ and $d = 50$. The first plot shows the one-dimensional ridge profile of \eqref{eq:ex_exactRidge} in black with the polynomial approximation computed using the generalized Gaussian rule from Algorithm \ref{alg:ridge_approx} denoted by the blue dashed line. The inverse CDF method using Gauss-Legendre quadrature from Tsilifis~\cite{Tsilifis2018} is shown in green for comparison. The absolute error of each method at each point is shown in the second plot, and the density $q(u)$ is in the third plot.  Visually, both methods appear to perform reasonably well in the center of the domain, where the induced density function is relatively large. However, by examining the absolute errors of each method, we see that Algorithm \ref{alg:ridge_approx} outperforms the pseudospectral approximation constructed using Gauss-Legendre quadrature by several orders of magnitude. Near the endpoints, the errors in each approximation method begin to increase, but $q(u)$ in these regions is many orders of magnitude smaller than in the middle of the domain.

\begin{figure*}[!ht]
\centering
\includegraphics[width=0.85\textwidth]{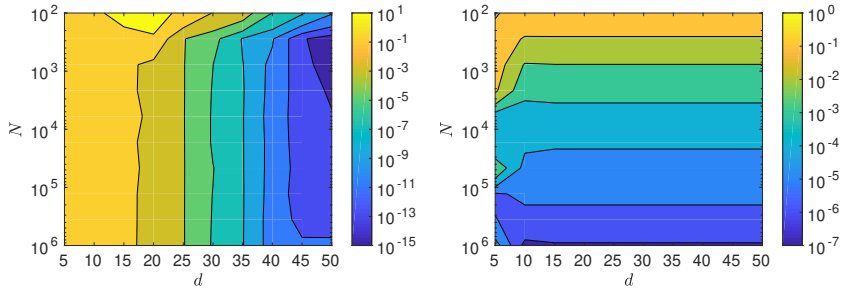}
\caption{Results from studies of the two levels of approximation in Algorithm \ref{alg:ridge_approx} applied to \eqref{eq:ex_exactRidge}. The left plot contains $L^2$ errors in the polynomial approximation for various values of $N$ and $d$. The right plot shows errors in the approximate integral of \eqref{eq:ex_exactRidge}.}
\label{fig:2d_convergence}
\end{figure*}

Next, we study the two levels of approximation in Algorithm \ref{alg:ridge_approx}: (i) the number $N$ of trapezoidal rule points used to construct the discrete approximation of $q(u)$ and (ii) the number $d$ of Lanczos iterations performed. The latter corresponds to the degree of the polynomial expansion as well as the number of Gaussian quadrature nodes. Figure \ref{fig:2d_convergence} shows the approximate $L^2$ norm of the error between $f$ and the pseudospectral approximations for varying values of $N$ and $d$. The error depends strongly on $d$. This is because a high-degree polynomial is required to fit the highly-oscillatory $f(\vx)$. The rightmost plot contains approximations of the integral of $f(\vx)$ using the first coefficient in the pseudospectral expansion. Here, we see a strong dependence on $N$. This is because integration errors in the Gaussian quadrature decay quickly with $d$. To improve the approximation, the discrete approximation of the density $q(u)$ must be improved.

\section{Extension to near-1D ridge functions}
\label{sec:near_ridge}

The approach presented in the previous section produces a univariate Gaussian quadrature for an exact ridge function. However, in many if not most applications, the exploratory analysis suggests that the computational model's input/output map can be merely approximated by a one-dimensional ridge function. Evidence for this near ridge structure appears in Figure \ref{fig:near_1d_ridge}, which shows the one-dimensional shadow plots for six different computational science applications. In each plot, the data shows relatively small deviation (relative to the range of the data) from a univariate function. 

Next we consider how Algorithm \ref{alg:ridge_approx} can be extended for the case when $f(\vx)$ is well approximated by a one-dimensional ridge function,
\begin{equation}
\label{eq:y_fx_approx_gaTx}
y \;=\; f(\vx) \;\approx\; g(u), 
\quad \text{where} \quad
u \;=\; \va^\top \vx .
\end{equation}
The best $L^2$ approximation of $f$ by $g$ is the expected value of the output conditioned on $u = \va^\top \vx$~\cite[Ch. 8]{Pinkus15}. That is, 
\begin{equation}
g(u) = \Exp{f(\vx) \,|\, \va^\top \vx = u} .
\end{equation}
By the tower property of conditional expectations~\cite{Billingsley86}, we can write the integral of $f$ as
\begin{equation}
\int f(\vx)\,p(\vx)\,\text{d}\vx \;=\; \int g(u)\,q(u)\,\text{d}u. 
\end{equation}
Thus, if we can evaluate the conditional expectation $g(u)$, then we can construct a univariate polynomial approximation and Gaussian quadrature for near ridge functions. We can measure the degree to which $f$ is a near one-dimensional ridge function using the conditional variance $\Var{f(\vx) \,|\, \va^\top \vx = u}$. This conditional variance is the error for the best one-dimensional ridge approximation of $f(\vx)$. 
%In shadow plots (see Figure \ref{fig:near_1d_ridge}), this variance is visually represented by the spread of the data for a fixed value of $u = \va^\top \vx$.

\begin{algorithm}[!ht] 
\caption{Univariate polynomial approximation of a near ridge function}
\label{alg:near_ridge_approx}
\textbf{Given:} function $f: \mathbb{R}^m \rightarrow \mathbb{R}$ and unit vector $\va \in \mathbb{R}^m \setminus \{ \mzero \}$ such that
$$y \;=\; f(\vx) \;\approx\; g(\va^\top \vx)$$
for some unknown $g: \mathbb{R} \rightarrow \mathbb{R}$ \\
\textbf{Assumptions:} $\vx \sim \sU ([-1, 1]^m)$ and $f$ is square-integrable with respect to the input density $p(\vx)$
\begin{enumerate}
\item Perform Steps 1-3 of Algorithm \ref{alg:ridge_approx}. 
\item For $j = 0,\dots,d$,
\begin{enumerate}
\item[(i)] Compute
\begin{equation}
\vxi_{j,0} \;=\; (1 - \gamma_j) \, \vx_\ell + \gamma_j \, \vx_r ,
\quad \text{for} \quad 
\gamma_j = (\lambda_j - u_\ell) \, / \, (u_r - u_\ell) ,
\end{equation}
where $\vx_\ell = \sign{-\va}$, $\vx_r = \sign{\va}$, $u_\ell = \va^\top \sign{-\va}$, and \\ $u_r = \va^\top \sign{\va}$.
\item[(ii)] Evaluate $f(\vxi_{j,0})$.
\item[(iii)] For $i = 1,\dots,M_j-1$,
\begin{enumerate}
\item[(a)] Choose a random vector $\vw \in \mathbb{R}^m \setminus \{ \mzero \}$ with $\vw^\top \va = 0$ and $||\vw||_2 = 1$.
\item[(b)] Choose a random value $t \in [-\sqrt{m}, \sqrt{m}]$.
\item[(c)] If $||\vxi_{j,i} + t \, \vw||_{\infty} <= 1$, then set $\vxi_{j,i+1} = \vxi_{j,i} + t \, \vw$. \\ Otherwise, repeat Step (b).
\item[(d)] Evaluate $f(\vxi_{j,i+1})$.
\end{enumerate}
\item[(iv)] Approximate the conditional expectation at the quadrature point,
\begin{equation}
g(\lambda_j) \;\approx\; \hat{g}(\lambda_j) \;=\; \frac{1}{M_j} \sum_{i=0}^{M_j-1} f(\vxi_{j,i}) 
\end{equation}
and the standard error $s_j \;=\; \hat{\sigma}_j/\sqrt{M_j}$, where $\hat{\sigma}_j$ is the standard deviation of the $f(\vxi_{j,i})$ for $i = 0,\dots,M_j-1$.
\end{enumerate}
\item Compute the pseudospectral coefficients 
\begin{equation}
\hat{g}_i \;=\; \sum_{j=0}^{d} \nu_j \, \hat{g}(\lambda_j) \, \phi_i (\lambda_j) .
\end{equation}
\item Define $\tilde{d}$ to be the largest index such that $|\hat{g}_i| < \sum_{j=0}^d s_j / (d+1)$ for all $i > \tilde{d}$.
\item Build the pseudospectral expansion 
\begin{equation}
y \;=\; g(u) \;\approx\; \sum_{i=0}^{\tilde{d}} \hat{g}_i \, \phi_i (u) .
\end{equation}
\end{enumerate}
\end{algorithm}

The difficulty in applying the methodology from section \ref{sec:exact_ridge} to \eqref{eq:y_fx_approx_gaTx} is how to compute $g(\lambda_j)$. Recall from \eqref{eq:1dQuad_to_mdQuad} that, in the case of the exact ridge function, we can transform $\lambda_j$ into a corresponding input $\vxi_j$ in the full-dimensional space and evaluate $f$ at this point, provided that $\vxi_j$ is in the domain of $f$. In the near ridge case, we want to approximate the average of $f(\vx)$ over all $\vx \in [-1, 1]^m$ such that $\lambda_j = \va^\top \vx$. We write the sample approximation of this conditional expectation as
\begin{equation}
\label{eq:approx_cond_exp}
g(\lambda_j) \;\approx\; \hat{g}(\lambda_j) \;=\; \frac{1}{M_j} \sum_{i=0}^{M_j-1} f(\vxi_{j,i}) ,
\end{equation}
where the $M_j$ input values $\vxi_{j,i} \in [-1,1]^m$ are sampled uniformly conditioned on $\va^\top \vxi_{j,i} = \lambda_j$. To compute \eqref{eq:approx_cond_exp}, we use a hit-and-run sampling algorithm~\cite{Lovasz1993}. Start by choosing a random (unit) direction $\vw \in \mathbb{R}^m \setminus \{ \mzero \}$ that is orthogonal to the ridge direction $\va$. We then pick a step size $t \in [-\sqrt{m}, \sqrt{m}]$, where this range is used to ensure that the step size covers the maximum possible range of the rotated $[-1, 1]^m$ hypercube. To obtain the $(i+1)$st conditional sample, we step from the previously drawn sample, $\vxi_{j,i+1} = \vxi_{j,i} + t \, \vw$, provided that $\vxi_{j,i+1} \in [-1,1]^m$. If this is not the case, then we choose a new random step size until a valid conditional sample is obtained.

The second issue in approximating a near-1D ridge function is that the sample approximation of the conditional expectation in \eqref{eq:approx_cond_exp} results in noisy estimates of $g(\lambda_j)$. For this reason, constructing an interpolating polynomial, as described in Section \ref{sec:exact_ridge}, is not the best approach given a restricted computational budget. We recommend truncating the pseudospectral polynomial expansion to fewer than $d$ terms (where $d$ is the number of Gaussian quadrature points) to avoid overfitting. For each Gaussian quadrature node $\lambda_j$, we have $\vxi_{j,i}$ for $i = 0,\dots,M_j-1$. Estimate the standard error in each approximation of the conditional expectation, $s_j \;=\; \hat{\sigma}_j/\sqrt{M_j}$, where $\hat{\sigma}_j$ is the standard deviation of the $f(\vxi_{j,i})$, $i = 0,\dots,M_j-1$. When $g(u)$ is smooth, we expect decay in the coefficients $g_i$ from \eqref{eq:1D_exp} for sufficiently large $i$~\cite{ATAP2013}. The pseudospectral coefficients, $\hat{g}_i$ from \eqref{eq:1D_pseudo_exp}, approximate the true coefficients. We suggest truncating the expansion at a degree $\tilde{d} \leq d$ polynomial, where $|\hat{g}_i| < \sum_{j=0}^d s_j / (d+1)$ for all $i > \tilde{d}$. This heuristic removes terms whose contribution to the expansion is smaller than the noise in the sample approximations $\hat{g} (\lambda_j)$. Algorithm \ref{alg:near_ridge_approx} formalizes the process for building a pseudospectral approximation of a near-1D ridge function.

\subsection{Error analysis}
\label{subsec:error_analysis}

We briefly explore the sources of error in the approximation and integration methods for near-1D ridge functions. For this discussion, let $p_d (u)$ denote the $d$-degree pseudospectral polynomial approximation of the ridge profile $g(u)$ constructed by Algorithm \ref{alg:near_ridge_approx}. Define the $m$-dimensional analogs of $p_d$ and $g$ by
\begin{equation}
p_{d,\vx} (\vx) \;=\; p_d(\va^\top \vx) 
\quad \text{and} \quad
g_{\vx} (\vx) \;=\; g(\va^\top \vx) .
\end{equation}
This enables us to study the $L^2$ errors without changing the behavior of the one-dimensional functions $p_d$ and $g$.

Using the triangle inequality, we can decompose the error in the one-dimensional approximation of $f$ as
\begin{equation} \label{eq:approx_error}
\left\| f - p_{d,\vx} \right\|_{L^2}
\;\leq\;
\left\| f - g_{\vx} \right\|_{L^2} +  \left\| g_{\vx} - p_{d,\vx}  \right\|_{L^2} .
\end{equation}
The two terms on the right-hand side of \eqref{eq:approx_error} represent the degree to which $f(\vx)$ is well approximated by a one-dimensional ridge function and the degree to which we approximate the ridge function by a polynomial, respectively. For exact ridge functions (see Section \ref{sec:exact_ridge}), $\left\| f - g_{\vx} \right\|_{L^2} = 0$ and the error only depends on the ability to approximate the ridge profile by a polynomial. For an accessible treatment of univariate polynomial approximation, see Trefethen's \emph{Approximation Theory and Approximation Practice}~\cite{ATAP2013}. For near-1D ridge functions, the ridge approximation error is non-zero. This implies that the error in the one-dimensional polynomial approximation of $f(\vx)$ cannot be reduced below this threshold. We show this behavior numerically in the next sections.

The integration error is tied to the approximation error. Using H{\"o}lder's inequality, 
\begin{equation}
\left\| f - g_{\vx} \right\|_{L^1}
\;\leq\;
\left\| f - g_{\vx} \right\|_{L^2} .
\end{equation}
Then, 
\begin{equation}
\left\| f - g_{\vx} \right\|_{L^1}
\;=\;
\int \left| f - g_{\vx} \right| \, p(\vx)\,\text{d}\vx
\;\geq\;
\left| \int (f - g_{\vx}) \, p(\vx)\,\text{d}\vx \right|
\;=\;
\left| \Exp{f} - \Exp{g_{\vx}} \right| .
\end{equation}
Thus,
\begin{equation}
\left| \Exp{f} - \Exp{g_{\vx}} \right|
\;\leq\;
\left\| f - g_{\vx} \right\|_{L^2} .
\end{equation}

\subsection{Numerical studies: Example 1}
\label{subsec:near_ridge_numerics}

In this section, we numerically study the behavior of Algorithm \ref{alg:near_ridge_approx} for approximating nearly one-dimensional ridge functions. We consider the function 
\begin{equation}
\label{eq:ex_nearRidge}
y \;=\; \sin \left( \frac{\pi}{5} \, (\va^\top \vx) \right) + \frac{1}{5} \, \cos \left( \frac{4 \, \pi}{5} (\va^\top \vx) \right) + \frac{1}{40} \vx^\top \mB \, \bf{1}
\qquad
\vx \in \mathbb{R}^{25} ,
\end{equation}
where $\va \in \mathbb{R}^{25} \setminus \{ \mzero \}$ defines the ridge-like structure, $\mB \in \mathbb{R}^{25 \times 24}$ contains an orthonormal basis for the subspace orthogonal to $\va$, and $\mathbf{1} \in \mathbb{R}^{24}$ is a vector of ones. We assume $\vx \sim \sU ([-1, 1]^{25})$.

\begin{figure*}[!ht]
\centering
\includegraphics[width=0.7\textwidth]{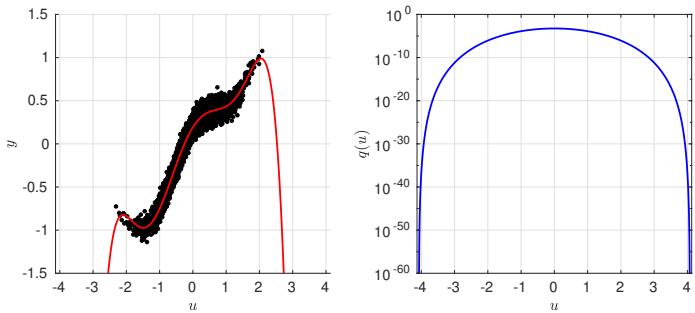}
\caption{The results of the extension to Algorithm \ref{alg:ridge_approx} applied to the approximate ridge function \eqref{eq:ex_nearRidge}. The left plot contains a show plot of $f(\vx)$ along the $u = \va^\top \vx$ axis with the $d = 11$ polynomial approximation on top of it. The right plot shows the density $q(u)$}
\label{fig:eyeball_norm2}
\end{figure*}

Figure \ref{fig:eyeball_norm2} shows the results of using the extension of Algorithm \ref{alg:ridge_approx} to \eqref{eq:ex_nearRidge}. The plot on the left is a shadow plot of evaluations of \eqref{eq:ex_nearRidge} against $u = \va^\top \vx$. The spread in the plot is due to variations in the 24 directions orthogonal to $\va$. The red line is the polynomial approximation of the ridge profile---$g(u) = \Exp{\,f(\vx) \,|\, \va^\top \vx=u}$. This is computed using $d = 11$ and $M = 50$ total function evaluations distributed among the 12 Gaussian quadrature nodes. The polynomial expansion is truncated at $\tilde{d} = 6$ to avoid overfitting the noise in the approximation of $g(\lambda_j)$. Note that fitting a polynomial of total degree 6 in 25 dimensions would require at least $M = {25+6 \choose 6} = 736,281$ to have a well-posed fitting problem.

\subsection{Numerical studies: Example 2}
\label{subsec:near_ridge_numerics2}

\begin{figure*}[!ht]
\centering
\includegraphics[width=0.5\textwidth]{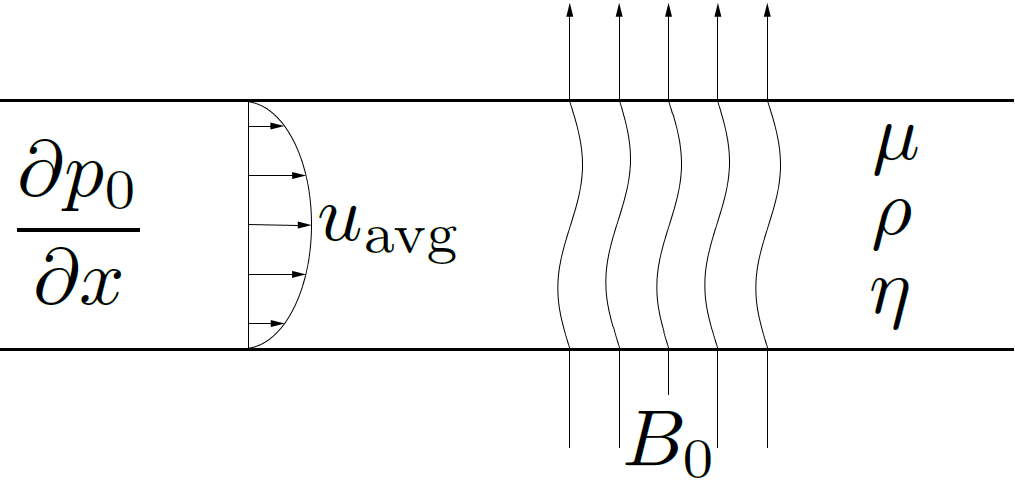}
\caption{An illustration of the Hartmann problem, which models the flow of an ionized fluid in the presence of a perpendicular magnetic field along an infinitely-long channel.}
\label{fig:Hartmann}
\end{figure*}

Next, we consider a physically-motivated problem in magnetohydrodynamics: the Hartmann problem~\cite{Cowling57}. The Hartmann problem models the flow of an ionized fluid in the presence of a perpendicular magnetic field along an  infinitely-long channel (see Figure \ref{fig:Hartmann}). The magnetic field acts as a resistive force on the flow while the fluid induces a horizontal component in the magnetic field. For this work, we consider the average flow velocity across the channel, denoted by $u_{\text{avg}}$, as the output of interest. This quantity can be written as
\begin{equation}
\label{eq:U_avg}
u_{\text{avg}} (\vx) = -\frac{\partial p_0}{\partial x} \frac{\eta}{B_0^2} \left( 1 - \frac{\ell B_0}{\sqrt{\eta \mu}} \, \text{coth} \left( \frac{\sqrt{\eta \mu}}{B_0 \ell} \right) \right) ,
\end{equation}
where the five input variables are described in Table \ref{tab:Hartmann_inputs}. The constant values $\ell$ and $\mu_0$ are the width of the channel and the magnetic permeability of free space (i.e., a universal constant), respectively. Recent work has shown that \eqref{eq:U_avg} exhibits approximate one-dimensional ridge structure with respect to the log-transformed inputs~\cite{Glaws17b}. Thus, we consider \eqref{eq:U_avg} with respect to the inputs
\begin{equation}
\vx \;=\; \bmat{\log(\mu) & \log(\rho) & \log(\partial p_0 / \partial x) & \log(\eta) & \log(B_0)}^\top .
\end{equation}
We assume a uniform density function over the range of values given in Table \ref{tab:Hartmann_inputs}.

\begin{table}[]
\centering
\caption{Descriptions and ranges of the five variable inputs to \eqref{eq:U_avg}. Note that the range of values is given in terms of the log of each input.}
\label{tab:Hartmann_inputs}
\begin{tabular}{ccc} \hline
Symbol & Description & Range of $\log (input)$ \\ \hline
$\mu$ & fluid viscosity & $[ \, \log(0.05), \, \log(0.2) \, ]$ \\
$\rho$ & fluid density & $[ \, \log(1), \, \log(5) \, ]$ \\
$\partial p_0 / \partial x$ & applied pressure gradient & $[ \, \log(0.5), \, \log(3) \, ]$ \\
$\eta$ & magnetic resistivity & $[ \, \log(0.5), \, \log(3) \, ]$ \\
$B_0$ & applied magnetic field & $[ \, \log(0.25), \, \log(1) \, ]$ \\ \hline
\end{tabular}
\end{table}

\begin{figure*}[!ht]
\centering
\includegraphics[width=1\textwidth]{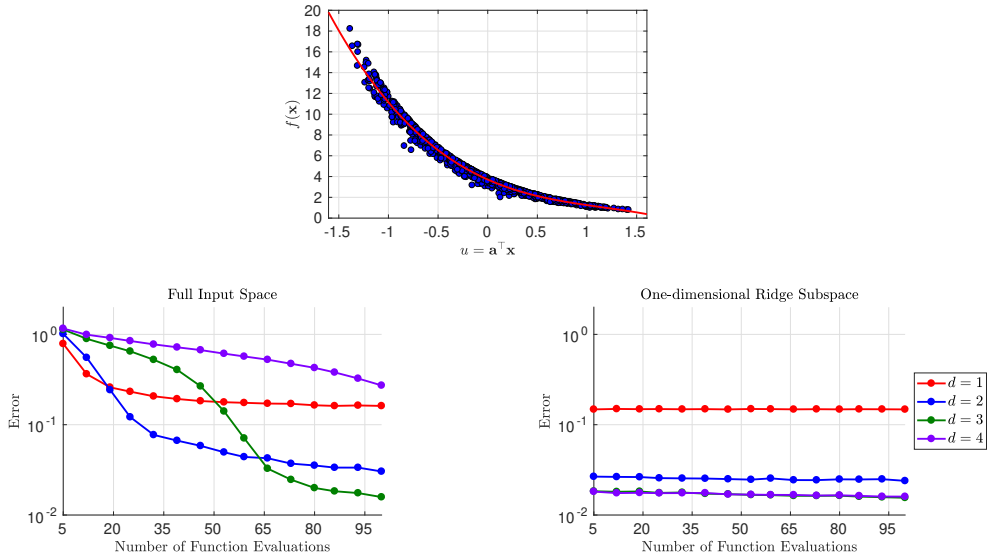}
\caption{The results of applying the extension of Algorithm \ref{alg:ridge_approx} to \eqref{eq:U_avg}. The top plot is a shadow plot of $u_{\text{avg}}$ against $u = \va^\top \vx$. The bottom plots show $L^2$ errors of polynomial approximations constructed on the full five-dimensional input space (on the left) and the one-dimensional ridge subspace (on the right).}
\label{fig:hartmann_nearRidge}
\end{figure*}

Figure \ref{fig:hartmann_nearRidge} contains the results of applying Algorithm \ref{alg:near_ridge_approx} to \eqref{eq:U_avg}. The top plot is a shadow plot of $u_{\text{avg}}$ against $u = \va^\top \vx$ overlaid with a polynomial approximation of the ridge profile. This approximation was constructed using $d = 4$ with a total computational budget of 100 function evaluations (20 for each of the 5 quadrature nodes). The bottom plots show the approximated relative $L^2$ errors in the polynomial approximations constructed on the full five-dimensional input space (on the left) and the one-dimensional ridge subspace (on the right). On the full input space, we use uniformly-sampled points from the $[-1,1]^5$ hypercube and construct the least-squares polynomial approximation using Legendre polynomials with $L^2$ regularization. These approximations perform poorly when the computational budget is limited and restrict our optimal choice of polynomial degree. This issue grows exponentially as the dimension of the given function increases. The one-dimensional approximation constructed using Algorithm \ref{alg:near_ridge_approx} achieves its optimal performance with very few function evaluations. For the restricted computational budget studied, the one-dimensional approximation outperforms its five-dimensional counterpart. Given a larger computational budget, the full polynomial approximation will be more accurate. This is because the one-dimensional approximation is limited by the accuracy of approximating $u_{\text{avg}} (\vx)$ by a ridge function (recall \eqref{eq:approx_error}). We estimate the degree to which $u_{\text{avg}} (\vx)$ is well approximated by a ridge function as
\begin{equation} \label{eq:uavg_near_ridge}
\left\| u_{\text{avg}} - g_{\vx} \right\|_{L^2}
\;\approx\;
1.29 \times 10^{-2} ,
\end{equation} 
which is approximately where the $L^2$ error in the one-dimensional ridge approximation begins to level out. The value of \eqref{eq:uavg_near_ridge} is approximated by taking 1000 uniformly-sampled points from $[-1,1]^5$ and transforming them into $u_i = \va^\top \vx_i$. At each $u_i$, we obtain 100 randomly-sampled points orthogonal to $\va$. Such a study is infeasible for an expensive computational model, but it is useful in explaining the behavior of the one-dimensional ridge approximations introduced in this paper.

\section{Summary and extensions}
\label{sec:conclusion}

We introduce a novel algorithm for estimating the integral and constructing a polynomial approximation of a one-dimensional ridge function based on Lanczos' method. In general, building a polynomial surrogate of an $m$-dimensional function suffers from the curse of dimensionality---an exponential increase in computational costs resulting from increases in $m$. We also introduce an approach to extending this algorithm to functions that are well-approximated by a one-dimensional ridge function.

We numerically study the new algorithm on several test problems, including exact and approximate ridge functions. We show that exploiting low-dimensional structure can result in exponential savings while maintaining accuracy. Additionally, we study the two-level approximation behavior of the algorithm: the first level is a discrete approximation of the induced density function $q(u)$, and the second level is the number of Lanczos iterations, which corresponds to the degree of the polynomial approximation of the ridge function as well as the number of Gaussian quadrature points. In studying the extension of the algorithm to nearly one-dimensional ridge functions, we show that we can quickly achieve the ridge approximation error using very few function evaluations.

In section \ref{sub:ridges}, we mentioned that projection pursuit regression~\cite{Friedman1981} uses a model for the conditional mean of the regression that is a sum of ridge functions, where each ridge profile in the sum is a smoothing spline. Classically, a maximum likelihood approach is used to estimate the directions and spline parameters. Once those parameters are estimated, our approach can be used to accurately estimate the integral of the projection pursuit regression model, since
\begin{equation}
\int \left(\sum_i g_i(\va_i^\top\vx)\right) \,p(\vx)\,\text{d}\vx \;=\; 
\sum_i \left(\int g_i(\va_i^\top\vx)\,p(\vx)\,\text{d}\vx\right).
\end{equation}
Essentially, one could repeat our process for each pair of $\va_i$ and $g_i$, and then add the contributions.

Finally, we mention the difficulty in extending our approach to more than one dimension. Some functions that arise in computational science may be well approximated by a \emph{generalized} ridge function $g(\mA^\top \vx)$, where $\mA\in\mathbb{R}^{m\times n}$ with $m>n$ and $g:\mathbb{R}^n\rightarrow\mathbb{R}$. Unfortunately, the vector $\mA^\top\vx$ does not in general contain independent components even when $\vx$'s components are independent. Therefore, straightforward tensor product extensions to univariate Gaussian quadrature are not appropriate. In principle, one could use linear programming extensions for Gaussian quadrature on convex domains~\cite{Ryu2015,Jakeman2018}. We leave such extensions for future work.

%\section*{Acknowledgments}
%THANKS!!!
\bibliographystyle{siamplain}
\bibliography{lanczos-1d-ridge}

\begin{thebibliography}{10}

\bibitem{Allen2004}
{\sc M.~Allen and K.~Maute}, {\em Reliability-based design optimization of
  aeroelastic structures}, Structural and Multidisciplinary Optimization, 27
  (2004), pp.~228--242, \url{https://doi.org/10.1007/s00158-004-0384-1}.

\bibitem{Beck2018}
{\sc J.~A. Beck, J.~M. Brown, A.~A. Kaszynski, and E.~B. Carper}, {\em Active
  subspace development of integrally bladed disk dynamic properties due to
  manufacturing variations}, in ASME Turbo Expo: Power for Land, Sea, and Air,
  Volume 7A: Structures and Dynamics, no.~51135, 2018, pp.~V07AT32A011--,
  \url{https://doi.org/10.1115/GT2018-76800}.

\bibitem{Berguin2014}
{\sc S.~H. Berguin and D.~N. Mavris}, {\em Dimensionality reduction using
  principal component analysis applied to the gradient}, AIAA Journal, 53
  (2014), pp.~1078--1090, \url{https://doi.org/10.2514/1.J053372}.

\bibitem{Billingsley86}
{\sc P.~Billingsley}, {\em Probability and Measure}, Wiley, New York, 1986,
  \url{https://www.wiley.com/en-us/Probability+and+Measure%2C+Anniversary+Edition-p-9781118341919}.

\bibitem{bungartz2004}
{\sc H.-J. Bungartz and M.~Griebel}, {\em Sparse grids}, Acta Numerica, 13
  (2004), pp.~147--269, \url{https://doi.org/10.1017/S0962492904000182}.

\bibitem{caflisch1998}
{\sc R.~E. Caflisch}, {\em {Monte Carlo and quasi-Monte Carlo methods}}, Acta
  Numerica, 7 (1998), pp.~1--49,
  \url{https://doi.org/10.1017/S0962492900002804}.

\bibitem{Cheney2000}
{\sc W.~Cheney and W.~Light}, {\em A Course in Approximation Theory}, American
  Mathematical Society, 2000.

\bibitem{Cohen12}
{\sc A.~Cohen, I.~Daubechies, R.~DeVore, G.~Kerkyacharian, and D.~Picard}, {\em
  Capturing ridge functions in high dimensions from point queries},
  Constructive Approximation, 35 (2012), pp.~225--243,
  \url{https://doi.org/10.1007/s00365-011-9147-6}.

\bibitem{Constantine15}
{\sc P.~G. Constantine}, {\em Active Subspaces: Emerging Ideas for Dimension
  Reduction in Parameter Studies}, SIAM, Philadelphia, 2015,
  \url{http://doi.org/10.1137/1.9781611973860}.

\bibitem{Constantine17a}
{\sc P.~G. Constantine and A.~Doostan}, {\em Time-dependent global sensitivity
  analysis with active subspaces for a lithium ion battery model}, Statistical
  Analysis and Data Mining, 10 (2017), pp.~243--262,
  \url{https://doi.org/10.1002/sam.11347}.

\bibitem{Constantine2014}
{\sc P.~G. Constantine, E.~Dow, and Q.~Wang}, {\em Active subspace methods in
  theory and practice: Applications to kriging surfaces}, SIAM Journal on
  Scientific Computing, 36 (2014), pp.~A1500--A1524,
  \url{https://doi.org/10.1137/130916138}.

\bibitem{Constantine17b}
{\sc P.~G. Constantine, A.~Eftekhari, J.~Hokanson, and R.~Ward}, {\em A
  near-stationary subspace for ridge approximation}, Computer Methods in
  Applied Mechanics and Engineering, 326 (2017), pp.~402--421,
  \url{https://doi.org/10.1016/j.cma.2017.07.038}.

\bibitem{constantine2015computing}
{\sc P.~G. Constantine, A.~Eftekhari, and M.~B. Wakin}, {\em {Computing active
  subspaces efficiently with gradient sketching}}, in IEEE 6th International
  Workshop on Computational Advances in Multi-Sensor Adaptive Processing
  (CAMSAP), Cancun, 2015, pp.~353--356,
  \url{https://doi.org/10.1109/CAMSAP.2015.7383809}.

\bibitem{Constantine12b}
{\sc P.~G. Constantine, M.~S. Eldred, and E.~T. Phipps}, {\em Sparse
  pseudospectral approximation method}, Computer Methods in Applied Mechanics
  and Engineering, 229 (2012), pp.~1--12,
  \url{https://doi.org/10.1016/j.cma.2012.03.019}.

\bibitem{Constantine15a}
{\sc P.~G. Constantine, M.~Emory, J.~Larsson, and G.~Iaccarino}, {\em
  Exploiting active subspaces to quantify uncertainty in the numerical
  simulation of the {HyShot II} scramjet}, Journal of Computational Physics,
  302 (2015), pp.~1--20, \url{https://doi.org/10.1016/j.jcp.2015.09.001}.

\bibitem{Constantine2015}
{\sc P.~G. Constantine and D.~F. Gleich}, {\em {Computing active subspaces with
  Monte Carlo}}, arXiv:1408.0545v2,  (2015),
  \url{https://arxiv.org/abs/1408.0545}.

\bibitem{Constantine2012}
{\sc P.~G. Constantine and E.~T. Phipps}, {\em A lanczos method for
  approximating composite functions}, Applied Mathematics and Computation, 218
  (2012), pp.~11751--11762, \url{https://doi.org/10.1016/j.amc.2012.05.009}.

\bibitem{Constantine15b}
{\sc P.~G. Constantine, B.~Zaharatos, and M.~Campanelli}, {\em Discovering an
  active subspace in a single-diode solar cell model}, Statistical Analysis and
  Data Mining, 8 (2015), pp.~264--273, \url{https://doi.org/10.1002/sam.11281}.

\bibitem{Cook94b}
{\sc R.~D. Cook}, {\em Using dimension-reduction subspaces to identify
  important inputs in models of physical systems}, in Proceedings of the
  Section on Physical and Engineering Sciences, American Statistical
  Association, Alexandria, VA, 1994, pp.~18--25,
  \url{http://users.stat.umn.edu/~rdcook/SDR/ASA94.pdf}.

\bibitem{Cook98}
{\sc R.~D. Cook}, {\em Regression Graphics: Ideas for Studying Regression
  through Graphics}, John Wiley \& Sons, Inc, New York, 1998,
  \url{http://doi.org/10.1002/9780470316931}.

\bibitem{Cortesi2017}
{\sc A.~Cortesi, P.~Constantine, T.~Magin, and P.~M. Congedo}, {\em Forward and
  backward uncertainty quantification with active subspaces: application to
  hypersonic flows around a cylinder}, hal-01592591,  (2017),
  \url{https://hal.inria.fr/hal-01592591/}.

\bibitem{Cowling57}
{\sc T.~G. Cowling and R.~B. Lindsay}, {\em Magnetohydrodynamics}, Physics
  Today, 10 (1957), p.~40, \url{http://doi.org/10.1063/1.3060498}.

\bibitem{Donoho00}
{\sc D.~L. Donoho}, {\em High-dimensional data analysis: The curses and
  blessings of dimensionality}, in AMS Conference on Math Challenges of the
  21st Century, 2000,
  \url{http://www-stat.stanford.edu/~donoho/Lectures/CBMS/Curses.pdf}.

\bibitem{eftekhari2017learning}
{\sc A.~Eftekhari, M.~B. Wakin, P.~Li, P.~G. Constantine, and R.~A. Ward}, {\em
  {Learning the second-moment matrix of a smooth function from point samples}},
  in 51st Asilomar Conference on Signals, Systems, and Computers, Asilomar, CA,
  2017, pp.~671--675, \url{https://doi.org/10.1109/ACSSC.2017.8335427}.

\bibitem{Fornaiser12}
{\sc M.~Fornasier, K.~Schnass, and J.~Vybiral}, {\em Learning functions of few
  arbitrary linear parameters in high dimensions}, Foundations of Computational
  Mathematics, 12 (2012), pp.~229--262,
  \url{http://doi.org/10.1007/s10208-012-9115-y}.

\bibitem{Friedman1981}
{\sc J.~H. Friedman and W.~Stuetzle}, {\em Projection pursuit regression},
  Journal of the American Statistical Association, 76 (1981), pp.~817--823,
  \url{https://doi.org/10.1080/01621459.1981.10477729}.

\bibitem{Gautschi04}
{\sc W.~Gautschi}, {\em Orthogonal Polynomials}, Oxford Press, Oxford, 2004,
  \url{https://global.oup.com/academic/product/orthogonal-polynomials-9780198506720}.

\bibitem{Ghanem1991}
{\sc R.~Ghanem and P.~Spanos}, {\em Stochastic Finite Elements: A Spectral
  Approach}, Springer-Verlag, New York, 1991,
  \url{https://www.springer.com/us/book/9781461277958}.

\bibitem{Gilbert16}
{\sc J.~M. Gilbert, J.~L. Jefferson, P.~G. Constantine, and R.~M. Maxwell},
  {\em Global spatial sensitivity of runoff to subsurface permeability using
  the active subspace method}, Advances in water resources, 92 (2016),
  pp.~30--42.

\bibitem{giles2015}
{\sc M.~B. Giles}, {\em {Multilevel Monte Carlo methods}}, Acta Numerica, 24
  (2015), pp.~259--328, \url{https://doi.org/10.1017/S096249291500001X}.

\bibitem{Glaws17a}
{\sc A.~Glaws, P.~G. Constantine, and R.~D. Cook}, {\em Inverse regression for
  ridge recovery: A data-driven approach for parameter space dimension
  reduction in computational science}, arXiv:1702.02227v1,  (2017),
  \url{https://arxiv.org/abs/1702.02227}.

\bibitem{Glaws17b}
{\sc A.~Glaws, P.~G. Constantine, J.~Shadid, and T.~M. Wildey}, {\em Dimension
  reduction in {MHD} power generation models: dimensional analysis and active
  subspaces}, Statistical Analysis and Data Mining, 10 (2017), pp.~312--325,
  \url{https://doi.org/10.1002/sam.11355}.

\bibitem{Golub10}
{\sc G.~H. Golub and G.~Meurant}, {\em Matrices, Moments, and Quadrature with
  Applications}, Princeton University, Princeton, 2010,
  \url{http://press.princeton.edu/titles/9104.html}.

\bibitem{Golub69}
{\sc G.~H. Golub and J.~H. Welsch}, {\em Calculation of {Gauss} quadrature
  rules}, Mathematics of Computation, 23 (1969), pp.~221--230,
  \url{https://doi.org/10.1090/S0025-5718-69-99647-1}.

\bibitem{Gramacy12}
{\sc R.~B. Gramacy and H.~Lian}, {\em Gaussian process single-index models as
  emulators for computer experiments}, Technometrics, 54 (2012), pp.~30--41.

\bibitem{grey2018active}
{\sc Z.~J. Grey and P.~G. Constantine}, {\em {Active subspaces of airfoil shape
  parameterizations}}, AIAA Journal, 56 (2018), pp.~2003--2017,
  \url{https://doi.org/10.2514/1.J056054}.

\bibitem{Griewank2008}
{\sc A.~Griewank and A.~Walther}, {\em Evaluating Derivatives: Principles and
  Techniques of Algorithmic Differentiation}, SIAM, Philadelphia, 2nd~ed.,
  2008, \url{https://doi.org/10.1137/1.9780898717761}.

\bibitem{Hastie2009}
{\sc T.~Hastie, R.~Tibshirani, and J.~Friedman}, {\em The Elements of
  Statistical Learning: Data Mining, Inference, and Prediction}, Springer, New
  York, 2nd~ed., 2009.

\bibitem{Higham2018}
{\sc C.~F. Higham and D.~J. Higham}, {\em Deep learning: An introduction for
  applied mathematicians}, arXiv:1801.05894v1,  (2018),
  \url{https://arxiv.org/abs/1801.05894}.

\bibitem{Hokanson18}
{\sc J.~M. Hokanson and P.~G. Constantine}, {\em Data-driven polynomial ridge
  approximation using variable projection}, SIAM Journal on Scientific
  Computing, 40 (2018), pp.~A1566--A1589,
  \url{https://doi.org/10.1137/17M1117690}.

\bibitem{Holodnak2018}
{\sc J.~Holodnak, I.~Ipsen, and R.~Smith}, {\em A probabilistic subspace bound
  with application to active subspaces}, SIAM Journal on Matrix Analysis and
  Applications, 39 (2018), pp.~1208--1220,
  \url{https://doi.org/10.1137/17M1141503}.

\bibitem{Hu2016}
{\sc X.~Hu, G.~T. Parks, X.~Chen, and P.~Seshadri}, {\em Discovering a
  one-dimensional active subspace to quantify multidisciplinary uncertainty in
  satellite system design}, Advances in Space Research, 57 (2016),
  pp.~1268--1279, \url{https://doi.org/10.1016/j.asr.2015.11.001}.

\bibitem{hyman2014}
{\sc J.~M. Hyman and M.~Tian}, {\em Accurate integration of high dimensional
  functions using polynomial detrending}, tech. report, Tulane University,
  2014.

\bibitem{Jakeman2018}
{\sc J.~D. Jakeman and A.~Narayan}, {\em Generation and application of
  multivariate polynomial quadrature rules}, Computer Methods in Applied
  Mechanics and Engineering, 338 (2018), pp.~134--161,
  \url{https://doi.org/10.1016/j.cma.2018.04.009}.

\bibitem{Jefferson16}
{\sc J.~L. Jefferson, J.~M. Gilbert, P.~G. Constantine, and R.~M. Maxwell},
  {\em Active subspaces for sensitivity analysis and dimension reduction of an
  integrated hydrologic model}, Computers \& Geosciences, 83 (2016),
  pp.~127--138, \url{http://doi.org/10.1016/j.cageo.2015.07.001}.

\bibitem{King2018}
{\sc R.~King, J.~Quick, C.~Adcock, and K.~Dykes}, {\em Active subspaces for
  wind plant surrogate modeling}, in 2018 Wind Energy Symposium. Kissimmee,
  Florida., American Institute of Aeronautics and Astronautics, 2018,
  \url{https://doi.org/10.2514/6.2018-2019}.

\bibitem{lanczos50}
{\sc C.~Lanczos}, {\em An iteration method for the solution of the eigenvalue
  problem of linear differential and integral operators}, Journal of Research
  of the National Bureau of Standards, 45 (1950), pp.~255--282.

\bibitem{LeMaitre2010}
{\sc O.~P. Le~Ma\^{i}tre and O.~M. Knio}, {\em Spectral Methods for Uncertainty
  Quantification}, Springer, New York, 2010,
  \url{https://doi.org/10.1007/978-90-481-3520-2}.

\bibitem{Lewis2016}
{\sc A.~Lewis, R.~Smith, and B.~Williams}, {\em Gradient free active subspace
  construction using morris screening elementary effects}, Computers \&
  Mathematics with Applications, 72 (2016), pp.~1603--1615,
  \url{https://doi.org/10.1016/j.camwa.2016.07.022}.

\bibitem{Li2018}
{\sc B.~Li}, {\em Sufficient Dimension Reduction: Methods and Applications in
  R}, CRC Press, Boca Raton, 2018.

\bibitem{Li2016}
{\sc W.~Li, G.~Lin, and B.~Li}, {\em Inverse regression-based uncertainty
  quantification algorithms for high-dimensional models: Theory and practice},
  Journal of Computational Physics, 321 (2016), pp.~259--278,
  \url{https://doi.org/10.1016/j.jcp.2016.05.040}.

\bibitem{Liesen13}
{\sc J.~Liesen and Z.~Strakos}, {\em Krylov Subspace Methods: Principles and
  Analysis}, Oxford Press, Oxford, 2013,
  \url{https://global.oup.com/academic/product/krylov-subspace-methods-9780199655410}.

\bibitem{Liu2017}
{\sc X.~Liu and S.~Guillas}, {\em Dimension reduction for {Gaussian} process
  emulation: An application to the influence of bathymetry on tsunami heights},
  SIAM/ASA Journal on Uncertainty Quantification, 5 (2017), pp.~787--812,
  \url{https://doi.org/10.1137/16M1090648}.

\bibitem{Lovasz1993}
{\sc L.~Lov\'{a}sz and M.~Simonovits}, {\em Random walks in a convex body and
  an improved volume algorithm}, Random Structures \& Algorithms, 4 (1993),
  pp.~359--412, \url{https://doi.org/10.1002/rsa.3240040402}.

\bibitem{lukaczyk2014active}
{\sc T.~W. Lukaczyk, P.~Constantine, F.~Palacios, and J.~J. Alonso}, {\em
  {Active subspaces for shape optimization}}, in 10th AIAA Multidisciplinary
  Design Optimization Conference, National Harbor, 2014, pp.~1--18,
  \url{https://doi.org/10.2514/6.2014-1171}.

\bibitem{Mayer15}
{\sc S.~Mayer, T.~Ullrich, and J.~Vyb\'iral}, {\em Entropy and sampling numbers
  of classes of ridge functions}, Constructive Approximation, 42 (2015),
  pp.~231--264, \url{http://doi.org/10.1007/s00365-014-9267-x}.

\bibitem{Novak1996}
{\sc E.~Novak and K.~Ritter}, {\em High dimensional integration of smooth
  functions over cubes}, Numerische Mathematik, 75 (1996), pp.~79--97,
  \url{https://doi.org/10.1007/s002110050231}.

\bibitem{Owen2013}
{\sc A.~B. Owen}, {\em Monte Carlo theory, methods and examples}, 2013,
  \url{http://statweb.stanford.edu/~owen/mc/}.

\bibitem{Pinkus15}
{\sc A.~Pinkus}, {\em Ridge Functions}, Cambridge University Press, 2015,
  \url{https://doi.org/10.1017/CBO9781316408124}.

\bibitem{Ryu2015}
{\sc E.~K. Ryu and S.~P. Boyd}, {\em Extensions of {Gauss} quadrature via
  linear programming}, Foundations of Computational Mathematics, 15 (2015),
  pp.~953--971, \url{https://doi.org/10.1007/s10208-014-9197-9}.

\bibitem{Seshadri2018}
{\sc P.~Seshadri, S.~Shahpar, P.~Constantine, G.~Parks, and M.~Adams}, {\em
  Turbomachinery active subspace performance maps}, Journal of Turbomachinery,
  140 (2018), pp.~041003--041003--11, \url{https://doi.org/10.1115/1.4038839}.

\bibitem{SmithUQ2013}
{\sc R.~C. Smith}, {\em Uncertainty Quantification: Theory, Implementation, and
  Applications}, SIAM, Philadelphia, 2013,
  \url{http://bookstore.siam.org/cs12/}.

\bibitem{Stieltjes84}
{\sc T.~J. Stieltjes}, {\em Quelques recherches sur la théorie des quadratures
  dites mécaniques}, Annales scientifiques de l'École Normale Supérieure, 1
  (1884), pp.~409--426, \url{http://eudml.org/doc/80911}.

\bibitem{SullivanUQ2015}
{\sc T.~J. Sullivan}, {\em Introduction to Uncertainty Quantification},
  Springer, New York, 2015, \url{https://doi.org/10.1007/978-3-319-23395-6}.

\bibitem{Traub98}
{\sc J.~F. Traub and A.~G. Werschulz}, {\em Complexity and Information},
  Cambridge University Press, Cambridge, 1998.

\bibitem{Trefethen2008}
{\sc L.~N. Trefethen}, {\em {Is Gauss quadrature better than
  Clenshaw-Curtis?}}, SIAM Review, 50 (2008), pp.~67--87,
  \url{https://doi.org/10.1137/060659831}.

\bibitem{ATAP2013}
{\sc L.~N. Trefethen}, {\em Approximation Theory and Approximation Practice},
  SIAM, Philadelphia, 2013.

\bibitem{Tripathy2016}
{\sc R.~Tripathy, I.~Bilionis, and M.~Gonzalez}, {\em Gaussian processes with
  built-in dimensionality reduction: Applications to high-dimensional
  uncertainty propagation}, Journal of Computational Physics, 321 (2016),
  pp.~191--223, \url{https://doi.org/10.1016/j.jcp.2016.05.039}.

\bibitem{Tsilifis2018}
{\sc P.~A. Tsilifis}, {\em Gradient-informed basis adaptation for {Legendre}
  chaos expansions}, Journal of Verification, Validation, and Uncertainty
  Quantification, 3 (2018), p.~011005, \url{https://doi.org/10.1115/1.4040802}.

\bibitem{Tyagi14}
{\sc H.~Tyagi and V.~Cevher}, {\em Learning non-parametric basis independent
  models from point queries via low-rank methods}, Applied and Computational
  Harmonic Analysis, 37 (2014), pp.~389--412,
  \url{https://doi.org/10.1016/j.acha.2014.01.002}.

\bibitem{Xiu02}
{\sc D.~Xiu and G.~E. Karniadakis}, {\em The {Wiener--Askey} polynomial chaos
  for stochastic differential equations}, SIAM Journal on Scientific Computing,
  24 (2002), pp.~619--644, \url{https://doi.org/10.1137/S1064827501387826}.

\end{thebibliography}
\end{document}